\documentclass[12pt]{article}
\usepackage{amsmath,amsfonts,amssymb}
\usepackage{mathrsfs}
\usepackage{dsfont}
\usepackage[utf8x]{inputenc}
\usepackage[T1]{fontenc}
\usepackage{lmodern} \usepackage{ucs}
\usepackage{fullpage}
\usepackage{float}
\usepackage{graphicx}
\usepackage{caption}

\usepackage{array}
\usepackage{multicol}
\usepackage{multirow}
\usepackage[format=hang, justification=raggedright]{caption}
\usepackage{subfig}

\usepackage{color}
\usepackage{authblk}
\usepackage{esint}
\usepackage{enumitem,hyperref}
\usepackage[normalem]{ulem}

\usepackage{tikz}
\usetikzlibrary{decorations}
\usetikzlibrary{decorations.pathreplacing}
\usetikzlibrary{decorations.pathmorphing}
\usepackage{tkz-tab}
\usetikzlibrary{shapes}
\tikzset{t style/.style={style=solid}}

\newtheorem{lemma}{Lemma}[section]
\newtheorem{theo}{Theorem}

\newtheorem{prop}[lemma]{Proposition}

\newtheorem{coro}[lemma]{Corollary}

\makeatletter
\def\namedlabel#1#2{\begingroup
    #2%
    \def\@currentlabel{#2}%
    \phantomsection\label{#1}\endgroup
}
\makeatother

\makeatletter
\renewcommand*{\eqref}[1]{%
  \hyperref[{#1}]{\textup{\tagform@{\ref*{#1}}}}%
}
\makeatother

\newcommand{\QED}{\mbox{}\hfill \raisebox{-0.2pt}{\rule{5.6pt}{6pt}\rule{0pt}{0pt}} \medskip\par}

\newcommand{\ds}{\displaystyle}
\newcommand{\ud}{\, {\mathrm{d}}}

\newcommand{\TumVar}{n}

\title {Fitting parameters of a Fokker-Planck-like equation with constraint}
\author[1]{Kevin Atsou\thanks{{\tt KokouKevin.Atsou@pfizer.com}}
}
\author[2]{Thierry~Goudon\thanks{{\tt thierry.goudon@univ-cotedazur.fr}
}}
\author[3]{Pierre-Emmanuel Jabin\thanks{{\tt pejabin@psu.edu}}
}
\affil[1]{\small Pfizer
}

\affil[2]{\small Universit\'e C\^ote d'Azur,   CNRS, LJAD
}

\affil[3]{\small Penn State University,  
Math. Dept.}

\date{}

\begin{document}
\maketitle
%\tableofcontents
\abstract{We 
analyse a Fokker-Planck like equation, driven by a scalar parameter in order to reach an integral constraint.
We exhibit criteria guaranteeing existence-uniqueness of a solution.
We also provide counter-examples.
This problem is motivated by an application to the immune control of tumor growth. 
}
%\sout{The model exhibits a possible control of the tumor growth by the  immune response.
%Nevertheless, the control is not complete in the sense that the asymptotic states keep residual tumors and activated immune cells.
%We investigate on numerical grounds the leading effects on the tumor-immune system interactions.
%In particular, we find that spatial hetererogeneities  
%of the distribution  of naive cells
%have a strong influence on the growth and control dynamics.}
%}

\vspace*{.5cm}
{\small
\noindent{\bf Keywords.}
Fokker-Planck equation. Constrained elliptic problems. Tumor growth.
 Equilibrium phase.%\\[.4cm]

\vspace*{.5cm}

\noindent{\bf Math.~Subject Classification.} 
35B99. %PDE, Qualitative properties of solutions, general
92C50, % Medical applications .
92C17, % Cell movement (chemotaxis, etc.) 
}

\section{Introduction}

We are interested in the following problem: 
given a certain ``confining'' (the meaning of which will be made precise later on) potential $\Phi:\mathbb R^N\to \mathbb R$, two positive constants
$\gamma, \ell>0$ 
and two non-negative
functions 
$\delta, S:\mathbb R^N\to \mathbb R$, we consider the PDE
\begin{equation}\label{def_FP}
\gamma u -\mu\nabla\cdot(\nabla \Phi u)-\Delta u=\mu S,\end{equation}
and we wonder whether the parameter $\mu>0$ can be  selected  so that the associated solution $u_\mu$ satisfies the constraint
 \begin{equation}\label{constraint}
 \ds\int_{\mathbb R^N} \delta u\ud x=\ell>0.\end{equation}
This question has been introduced 
in \cite{aabg_pub}, motivated by the modeling of the immune response to tumor growth, 
in order to explain equilibrium phases where the tumor is kept under control by the action of the immune cells.
Numerical simulations show that  
 the formation of equilibria, and thus  the existence and stability of solutions of \eqref{def_FP}-\eqref{constraint}, 
 is a quite robust phenomenon, see also \cite{aabg_fr}.
However, the analysis provided in \cite{aabg_pub}, by means of the implicit function theorem, is  
restricted 
to small values of the constraint parameter $\ell$, see Theorem~\ref{theo}.
We wish to extend the existence-uniqueness of the pair $(\mu, u_\mu)$ 
satisfying \eqref{def_FP}-\eqref{constraint}, associated to any $\ell\geq 0$.
In fact, our analysis shows that this requires further ``compatibility'' conditions between 
the potential $\Phi$, the source $S$ and the constraint function $\delta$.
We provide counter-examples explaining the  
role of these conditions.
Our arguments, which are likely of interest beyond the original application to tumor-immune system interactions, rely on properties of the underlying Fokker-Planck operator, moment propagation and duality reasonings.
To be more specific, the flavor of our result can be summarized as follows:
we assume the following confining assumption
\begin{equation}
\label{hPhi2}\tag{{\bf{H$\Phi$}}}
\begin{array}{l}
\textrm{$\Phi$ is $C^2$, $\nabla \Phi(0)=0$ and there exists  constants $\Lambda,M>0$ such that  }
\\
\textrm{the hessian matrix of $\Phi$ satisfies, for any $x\in \mathbb R^N$,
$\Lambda\mathbb I\leq D^2_{ij}\Phi(x)\leq M\mathbb I$,}
\end{array}
\end{equation}
and the following assumptions on the data
\[\begin{array}{l}
\textrm{$S:\mathbb R^N\to [0,\infty)$ is $C^\infty$ and compactly supported,}
\\
\textrm{$\delta: \mathbb R^N\to [0,\infty)$ is $C^1$, bounded  
and takes positive values on the neighborhood of $0$}.
\end{array}\]
Then, for any $\ell\geq 0$, we can find a $\mu\in [0,\infty)$ such that \eqref{def_FP}-\eqref{constraint} holds.
(Slightly more general assumptions on $S$ will considered below, see Theorem~\ref{C:ex}.)
Uniqueness is established in the \emph{radially symmetric} framework, see Theorem~\ref{mainth}, assuming
 the convexity/monotonicity properties
\begin{equation}\label{cx}\tag{{\bf{HCM}}}
\partial_r \Phi\geq 0,\qquad \partial^2_r \Phi\geq 0,\qquad \partial_r\delta\leq 0.\end{equation}

The paper is organized as follows.
First, Section~\ref{Sec:Motiv}  motivates the problem \eqref{def_FP}-\eqref{constraint} by rapidly coming back to the modeling 
introduced in \cite{aabg_pub}.
Next, Section~\ref{Sec:FP} discusses some properties of the Fokker-Planck operator $\nabla\cdot(\nabla\Phi u)+\Delta u$
which arises in \eqref{def_FP} and are crucial for the analysis.
In Section~\ref{Sec:Asymp}, we study in details the behavior of the constraint functional $$\mathscr F:\mu\longmapsto \int_{\mathbb R^N} \delta u_\mu\ud x$$ for small and large $\mu$'s.
Finally, in Section~\ref{Sec:Monot} we analyze the monotonicity of the mapping $\mathscr F$, depending on assumptions on the data $  
\Phi, S, \delta$. The analytical results are further illustrated 
by a few numerical examples.

\section{Motivation}
\label{Sec:Motiv}

We remind the reader the modeling 
principles 
that lead to the problem \eqref{def_FP}-\eqref{constraint}.
The earliest stages of tumor growth can be described 
through the evolution of the density of tumor cells
$(t,z) \mapsto \TumVar(t,z)$:  the integral   $\int_a^b z \TumVar(t,z)\ud z$ gives 
 the volume of the tumor occupied at time $t$ by tumor cells  having their size $z$ in the interval $(a,b)$.
 It is governed by two phenomena:
 a natural growth, embodied in the rate $z\mapsto V(z)\geq 0$ 
 and cell division mechanisms, where a cell with size $z'$ divides into cells with respective sizes $  z$ and $z'-z$.
 The latter depend on  
 the  frequency  of division   $z\mapsto a(z)$ 
and the  size-distribution $k(z\vert z')$ from  the division of a tumor cell with  size $z'$.    
Therefore, 
without any further interaction, the evolution of the tumor cells obey the initial-boundary value problem:
\begin{equation}\label{evol_n}
\partial_t n+\partial_z (Vn)=Q(n), \qquad n(0,z)=n_0(z),\qquad n(t,0)=0,\end{equation}
with 
\[
  Q(n)(t,z) = -a(z)n(t,z) + \int_z^{\infty}a(z')k(z\vert z')n(t,z')\ud z '.
\]
A basic example of such cell- division operator is given by the 
binary division operator
\[
    Q(n)(t,z) = 4a(2z)n(t,2z) -a(z)n(t,z).
\]
In any cases, the assumption on the kernel $k$ are such that 
\begin{itemize}
\item The total number of tumor cells is non decreasing
\[
    \dfrac{\ud}{\ud t}\ds\int_0^{\infty} n(t,z)\ud z   \geq 0,
\]
        \item The total mass of the tumor is non decreasing
\[
    \dfrac{\ud}{\ud t} \ds\int_0^{\infty} z n(t,z)\ud z  = \ds\int_0^\infty V(z) n(t,z)\ud z \geq 0.
\]
  \end{itemize}
  Note that the former is due to cell division, the latter to the natural growth.

  A remarkable fact about this growth-division equation is the existence of an eigenpair $(\lambda, N)$,
with  $\lambda >0$ and 
$z\geq 0 \mapsto N(z)$ taking non-negative values,
that satisfy
\begin{equation}
\label{tumorgrowthEigenProblem}
\left\{
\begin{array}{l}
 \partial_z (V N)   -Q( N) +  \lambda  N = 0 \text{ for  } z \geq  0 \\ 
 N(0) = 0, \qquad  N(z)>0 \ \text{ for $z>0$},\qquad  \ds\int_{0}^{+\infty}  N(z)\ud z  = 1.
\end{array}\right.\end{equation}
We refer the reader to \cite{DoGa, Mich1,PerBk} 
for precise assumptions and statements with proofs relying on suitable applications of the Krein-Rutman theorem.
Note that, in the specific case where $a, V$ are constant and $Q$ is the binary division operator, we have $\lambda=a$
and the profile $N$ is explicitly known, \cite{Bac, PerBk,PeRy}.
Dedicated numerical methods to compute 
the eigenpair are presented in 
\cite{aabg_fr}.
Furthermore, it can be shown that this eigenstate drives the large time behavior
of the Cauchy problem for \eqref{evol_n}: we have 
$
n(t,z) \sim_{t\to \infty} \nu_0 e^{\lambda t} \overline N(z)
$
where $\nu_0$ is a constant determined by the initial condition, see \cite{DGL,Mich1,MMP}.

The modeling of immune response adopted in \cite{aabg_pub} assumes that the displacement of the immune cells
holds at a larger scale, described by a space variable $x\in \mathbb R^N$, 
while the tumor 
is attached at a given location $x_0$.
% (or, if there are several 
%tumors, each deseases centers is discribed by such an equation).
The  immune
  cells are activated from a reservoir of resting cells, 
  and their motion is driven by diffusion 
  and chemotaxis directed towards the tumor.
  The strength of both the activation and the directed drift 
  depends on the total tumor mass
  \[\mu(t)=\ds\int_0^\infty z n(t,z)\ud z.\]
    Let $\Phi:\mathbb R^N\rightarrow \mathbb R$
    be a potential, intended to create 
    an attractive force towards the tumor location $x_0$ (and from now on, wlog, we suppose $x_0=0$).
    The time evolution of the concentration  of immune cells
    $C:(0,\infty)\times \mathbb R^N\rightarrow [0,\infty)$  is governed by 
    \[\partial_ t C-\Delta C-\mu\nabla\cdot(\nabla\Phi C) =\mu S-\gamma C,
    \]
    where $S: \mathbb R^N\rightarrow[0,\infty)$ describes the source of resting immune cells, and $\gamma>0$ is the 
    natural death rate 
    of the immune cells.
    The action of the immune cells on the tumor cells is taken into account through a death term 
    \[-\ds\int_{\mathbb R^N} \delta(x)C(t,x)\ud x,\]
    in the right hand side of \eqref{evol_n}
    where $\delta: \mathbb R^N\rightarrow[0,\infty)$ is intended 
    to describe the 
    killing effects on the tumor cells; it acts as a mollified delta-Dirac at $x_0=0$.
    
    Performing simulations of the coupled problem, one observes the formation of an equilibrium phase, 
    with a residual tumor, having a positive mass, controlled by the action of the immune cells.
    Such an equilibrium can be explained by coming back to the eigenproblem \eqref{tumorgrowthEigenProblem}:
       the death term induced by the immune cells is expected to  counterbalance the natural growth rate of the cell-division equation. 
       The other way around, we expect that $C(t,x)$ tends to $u(x) $, solution of \eqref{def_FP}-\eqref{constraint}
        as $t$ goes to $\infty$, with $\ell=\lambda$, the eigenvalue determined 
       by \eqref{tumorgrowthEigenProblem}.
       We refer the reader to \cite{aabg_pub, aabg_fr} for numerical illustration of such  a behavior, which seems very robust.
       Moreover, 
   this interpretation of the equilibrium phase by means of an eigenvalue problem permits to compute a priori 
   the final mass of the tumor, the parameter $\mu$ in \eqref{def_FP}-\eqref{constraint},  given the biological parameters \cite{aabg_fr}.
 Therefore, we wonder whether we can find a solution
   of the constrained problem \eqref{def_FP}-\eqref{constraint},
   for any value of the constraint parameter $\ell$, since we wish that the latter coincides with 
   the eigenvalue $\lambda$ in  \eqref{tumorgrowthEigenProblem}, determined by the parameters
   of the tumor growth equation.
   Unfortunately, a direct reasoning justifies this interpretation only
   for small values of $\ell$ (which 
   thus means a small eigenvalue $\lambda$ in  \eqref{tumorgrowthEigenProblem} and
   in terms of modelling means a tumor with low aggressiveness).

\begin{theo}\label{theo}
%Let $x\mapsto p S(x)\in L^2(\Omega)$ be  a non-negative function.
%Let $\Phi$ be the solution of $$-\nabla\cdot(\mathcal K\nabla \Phi)=\sigma - \ds\frac{1}{|\Omega|}\ds\int_\Omega \sigma(y)\ud y
%,$$
%endowed with the homogeneous  Neumann  boundary condition.
 If $\ell>0$ is small enough, there exists a unique $\mu(\ell)>0$ such that $u_{\mu(\ell)} $, solution of 
the stationary equation \eqref{def_FP},
satisfies \eqref{constraint}.
\end{theo}

\noindent
{\bf Proof.} The framework slightly differs from \cite{aabg_pub} which deals with the problem set in a bounded domain, endowed with appropriate boundary conditions. The argument uses the results in Proposition~\ref{prop:LxM} and Lemma~\ref{lem:C1}, detailed  below.
We are searching for  the zeroes of the mapping 
\[\mathscr X:(\ell,\mu)\in [0,\infty)\times [0,\infty)\longmapsto \ds\int _{\mathbb R^N} \delta u_{\mu} \ud x-\ell\]
where $u_{\mu} $ is the solution of \eqref{def_FP} associated to $\mu$ and knowing that
$\mathscr X(0,0)=0$, since $u_{0}=0$.
We have $\partial_{\mu} \mathscr X(\ell,\mu)=\int _\Omega \delta u'_{\mu} \ud x$, with $u'_{\mu}$, solution of 
\[
\gamma u'-\Delta  u'-\mu\nabla\cdot (u'\nabla\Phi)=
S
+\nabla\cdot (u_{\mu}\nabla\Phi)
\]
(see  Lemma~\ref{lem:C1} below).
Since $u_0=0$ and $S\geq 0$, we get $u'_{0}>0$ (see Proposition~\ref{prop:LxM} below).
It follows that $\partial_{\mu} \mathscr X(0,0)=\int _{\mathbb R^N} \delta u'_{0} \ud x>0$.
The implicit function theorem tells us that there exists $\ell
_\star>0$ 
and a mapping $\mu:\ell\in [0,\ell_\star)\mapsto\mu(\ell)$ such that $\mathscr X(\ell,\mu(\ell))=0$ holds
for any $\ell\in [0,\ell_\star)$, which means that $u_\ell$ satisfies \eqref{def_FP}-\eqref{constraint}.
Observe that
$$
\begin{array}{l}
\partial_{\ell} \mathscr X(\ell,\mu(\ell))+ \mu'(\ell)\partial_{\mu} \mathscr X(\ell,\mu(\ell))
=-1 +\mu'(\ell)\partial_{\mu} \mathscr X(\ell,\mu(\ell))=0\end{array}$$ holds
with $\partial_{\mu} \mathscr X(0,0)>0$. 
Hence, $\ell\mapsto \mu(\ell)$ is increasing on the neighborood of $\ell=0$, and it thus takes positive values.

Note that the argument cannot be extended for any $\ell$, since we do not have a direct knowledge on the sign of 
$\nabla\cdot (u_{\mu}\nabla\Phi)$ for $\mu\neq 0$. However, the proof does not use the confining feature of the potential $\Phi$.
Hence, we are going to develop a viewpoint that further exploits these properties.
\QED

\section{Fundamental properties of the operator $L_\mu$}
\label{Sec:FP}

Let us make the following assumption on the potential $\Phi$:
\begin{equation}
\label{hPhi1}
\begin{array}{l}
\textrm{for any $x\in \mathbb R^N$, we have $\Phi(x)\geq 0$},
\\
\textrm{ and, for any $\mu>0$,  $x\mapsto M_\mu(x)=e^{-\mu \Phi(x)}\in L^1(\mathbb R^N)$}.
\end{array}
\end{equation}
As a  matter of fact, the latter integrability property is guaranteed by the %following
 strengthened convexity condition \eqref{hPhi2}.
%assuming $\Phi\in C^2$,  
% \begin{equation}
%\label{hPhi2}
%\begin{array}{l}
%\textrm{$\nabla \Phi(0)=0$ and there exists a constant $\Lambda>0$ such that  }
%\\
%\textrm{the hessian matrix of $\Phi$ satisfies, for any $x\in \mathbb R^N$,
%$D^2_{ij}\Phi(x)\geq \Lambda\mathbb I$.}
%\end{array}
%\end{equation}
These conditions describe the confining feature of the potential, having an attractive effect towards $x=0$, which  is a strict global minimizor of the potential.
Then, we introduce the Fokker-Planck operator
\begin{equation}
\label{def_fp}
L_\mu u =\mu\nabla\cdot(\nabla\Phi u)+\Delta u
\end{equation}
and its 
adjoint operator (defined with the standard $L^2$ inner product)
\begin{equation}
\label{def_fpst}
L_\mu^\star \psi =-\mu\nabla\Phi \cdot\nabla \psi+\Delta \psi.
\end{equation}
It is convenient to recast these operators by making the function $M_\mu$ appear
\[
L_\mu u=\nabla\cdot\left(M_\mu \nabla\left(\ds\frac u{M_\mu}\right)\right),
\qquad 
L_\mu^\star \psi=\ds\frac{1}{M_\mu}\nabla\cdot\left( M_\mu \nabla\psi \right).
\]
Accordingly, we observe that
\begin{equation}\label{vari_L}
-\ds\int_{\mathbb R^N}  \ds\frac u{M_\mu}L_\mu u\ud x
=\ds\int_{\mathbb R^N}   M_\mu \left|\nabla\left(\ds\frac u{M_\mu}\right) \right|^2\ud x
\geq 0.
\end{equation}
We guess from this relation that the kernel of $L_\mu$ is spanned by $M_\mu$; this is indeed the case, 
as a consequence of the following Sobolev inequality,  owing to \eqref{hPhi2}, 
\begin{equation}\label{coer_L}\begin{array}{lll}
\ds\int_{\mathbb R^N}   M_\mu \left|\nabla\left(\ds\frac u{M_\mu}\right) \right|^2 \langle M_\mu\rangle\ud x&\geq& 
2\Lambda\mu \ds\int_{\mathbb R^N}   \left| u-\langle u\rangle  \ds\frac{M_\mu}{\langle M_\mu\rangle} \right|^2\ds\frac{\langle M_\mu\rangle\ud x}{M_\mu}
\\&\geq& 2\Lambda\mu\left(  \ds\int_{\mathbb R^N}   \left| u-\langle u\rangle \ds\frac{M_\mu}{\langle M_\mu\rangle} \right|\ud x\right)^2
,
\end{array}\end{equation}
holds, where $\langle u\rangle=\int_{\mathbb R^N} u\ud x$, see  \cite[condition (A2),  Corollary 2.18]{AMTU}.
Similarly, we have
\[
-\ds\int_{\mathbb R^N} M_\mu \psi L^\star_\mu \psi\ud x
=\ds\int_{\mathbb R^N} M_\mu \left|  \nabla\psi \right|^2\ud x
\geq 0.
\]
\begin{prop} 
\label{prop:LxM}
The following assertions hold:
\begin{itemize}
\item[i)]  $\mathrm{Ker}(L_\mu)=\mathrm{Span}(M_\mu)$ and $\mathrm{Ker}(L^\star_\mu)=\mathrm{Span}(\mathbf 1)$.
\item[ii)]  Let $\gamma>0$. For any $S\in L^2(\mathbb R^N,\frac{\ud x}{M_\mu})$, 
there exists a unique solution $u\in L^2(\mathbb R^N,\frac{\ud x}{M_\mu})$, with 
$\nabla\frac u{M_\mu}\in L^2(\mathbb R^N,M_\mu \ud x)$, of 
$(\gamma\mathbb I-L_\mu)u=S$.
Moreover, if $S\geq 0$, then $u\geq 0$.
\item[iii)]  Let $\gamma>0$. For any $\delta \in L^2(\mathbb R^N,M_\mu \ud x)$, 
there exists a unique solution $\psi\in L^2(\mathbb R^N,M_\mu \ud x)$, with 
$\nabla\psi\in L^2(\mathbb R^N,M_\mu \ud x)$ of 
$(\gamma\mathbb I-L^\star_\mu)\psi=\delta$.
Moreover, if $\delta \geq 0$, then $\psi\geq 0$.
\end{itemize}
\end{prop}

\noindent
{\bf Proof.}
The first item is a direct consequence of \eqref{vari_L} and \eqref{coer_L}.
Next, we simply apply the Lax-Milgram theorem (or, in the present context the Riesz theorem) in the Hilbert space
\[
H=\Big\{u\in  L^2\Big(\mathbb R^N,\frac{\ud x}{M_\mu}\Big),\ \nabla\Big( \ds\frac u{M_\mu}\Big)\in L^2(\mathbb R^N,M_\mu \ud x)\Big\}\]
to solve the variational problem: to find $u\in H$, such that, for any $v\in H$, we have
\[
\gamma \ds\int_{\mathbb R^N}uv \ds\frac{\ud x}{M_\mu}+\ds\int_{\mathbb R^N}   M_\mu \nabla\left(\ds\frac u{M_\mu}\right) 
\cdot \nabla\left(\ds\frac v{M_\mu}\right)
\ud x=\ds\int_{\mathbb R^N} Sv\ds\frac{\ud x}{M_\mu}.\]
We obtain the sign property by using $v=u_-=\min(0,u)$ as trial function in the variational formulation: it yields
\[
\gamma \ds\int_{\mathbb R^N} u^2_- \ds\frac{\ud x}{M_\mu}+\ds\int_{\mathbb R^N}   M_\mu \left|\nabla\left(\ds\frac {u_-}{M_\mu}\right) \right|^2\ud x=\ds\int_{\mathbb R^N} Su_-\ds\frac{\ud x}{M_\mu}\leq 0\]
when $S$ takes non-negative values. It implies $u_-=0$ a. e.  
  A similar argument applies readily to the adjoint problem.
Note that the variational formulation provides the estimate, for $u$ solution of $(\gamma-L_\mu)u=S$, 
  \[
  \gamma\ds\int_{\mathbb R^N}\ds\frac{|u|^2}{M_\mu}\ud x
  +\ds\int_{\mathbb R^N}   M_\mu \left|\nabla\left(\ds\frac {u}{M_\mu}\right) \right|^2\ud x=\ds\int_{\mathbb R^N} Su\ds\frac{\ud x}{M_\mu}
  \leq  \ds\frac\gamma2 \ds\int_{\mathbb R^N}\ds\frac{|u|^2}{M_\mu}\ud x+\ds\frac1{2\gamma} \ds\int_{\mathbb R^N}\ds\frac{|S|^2}{M_\mu}\ud x\]
  hence
  \[\ds\int_{\mathbb R^N}\ds\frac{|u|^2}{M_\mu}\ud x
  +\ds\frac2\gamma\ds\int_{\mathbb R^N}   M_\mu \left|\nabla\left(\ds\frac {u}{M_\mu}\right) \right|^2\ud x\leq 
  \ds\frac1{\gamma^2} \ds\int_{\mathbb R^N}\ds\frac{|S|^2}{M_\mu}\ud x\]
  that will be repeatedly used in what follows.
\QED

\section{Asymptotic behavior of $ \mu\mapsto \mathscr F(\mu)$} 
\label{Sec:Asymp}

As a warm up, we start by checking that $\mathscr F$ is well-defined.

\begin{lemma}
Suppose \eqref{hPhi2} holds.
Let $\delta:  \mathbb R^N\to [0,\infty)$ be a bounded function and let $S:\mathbb R^N\to [0,\infty)$ be such that 
$S\in L^2(\mathbb R^N, 
\frac{\ud x}{M_\mu})$ for any $\mu\geq 0$.
The function $\mathscr F$ takes value in $[0,\infty)$, and is continuous on $[0,\infty)$.
\end{lemma}

\noindent
{\bf Proof.}
We denote by $u_\mu$ the solution of \eqref{def_FP}, obtained by applying Proposition~\ref{prop:LxM}.
Since $u_0=0$, we have $\mathscr F(0)=0$, and, by virtue of Proposition~\ref{prop:LxM}-ii),  for any $\mu\geq  0$, 
$u_\mu\geq 0$, so that $\mathscr F(\mu)\geq 0$. Moreover, for $\mu>0$, we have
\[\begin{array}{lll}
\mathscr F(\mu)&\leq& \|\delta\|_{L^\infty(\mathbb R^N)}\ds\int_{\mathbb R^N} u_\mu\ud x\leq 
\|\delta\|_{L^\infty(\mathbb R^N)}
\|M_\mu\|_{L^1(\mathbb R^N)}^{1/2}
\left(\ds\int_{\mathbb R^N} \ds\frac{|u_\mu|^2}{M_\mu}\ud x\right)^{1/2}
\\&\leq& 
\mu\ds\frac{\|\delta\|_{L^\infty(\mathbb R^N)}}{\gamma^2}\|M_\mu\|_{L^1(\mathbb R^N)}^{1/2}
\left(\ds\int_{\mathbb R^N} \ds\frac{|S|^2}{M_\mu}\ud x\right)^{1/2}<\infty.\end{array}\]
By integrating the equation $(\gamma -L_\mu)u_\mu=\mu S$, we get
\begin{equation}\label{massmu}
\gamma\ds\int_{\mathbb R^N} u_\mu\ud x=\mu \ds\int_{\mathbb R^N} S\ud x.\end{equation}
It implies  the continuity of $\mathscr F$ at $\mu=0$: $\lim_{\mu\to 0}\mathscr F(\mu)=0$.
Next, let $0< \mu_1<\mu_2<\mu_*<\infty$ and $\epsilon=u_{\mu_2}-u_{\mu_1}$. It satisfies
\[
\gamma \epsilon -L_{\mu_1}\epsilon=(\mu_2-\mu_1)S+ (\mu_2-\mu_1)\nabla\cdot(\nabla\Phi u_{\mu_2}).\]
We deduce that 
\[\begin{array}{l}
\gamma\ds\int_{\mathbb R^N} \ds\frac{\epsilon^2}{M_{\mu_1}}\ud x
+\ds\int_{\mathbb R^N} M_{\mu_1}\Big|\nabla\Big(\ds\frac{\epsilon^2}{M_{\mu_1}}\Big)\Big|^2\ud x
\\
=
(\mu_2-\mu_1)\ds\int_{\mathbb R^N}S\ds\frac{\epsilon}{M_{\mu_1}}\ud x- 
(\mu_2-\mu_1)\ds\int_{\mathbb R^N}\nabla\Phi u_{\mu_2} \nabla\Big(\ds\frac\epsilon{M_{\mu_1}}\Big)\ud x
\\
\leq \ds\frac\gamma 2 \ds\int_{\mathbb R^N} \ds\frac{\epsilon^2}{M_{\mu_1}}\ud x
+ \ds\frac{(\mu_2-\mu_1)^2}{2\gamma}\ds\int_{\mathbb R^N}\ds\frac{S^2}{M_{\mu_1}}\ud x
\\\qquad+ \ds\frac12\ds\int_{\mathbb R^N} M_{\mu_1}\Big|\nabla\Big(\ds\frac{\epsilon^2}{M_{\mu_1}}\Big)\Big|^2\ud x
+\ds\frac{(\mu_2-\mu_1)^2}{2}\ds\int_{\mathbb R^N}|\nabla\Phi|^2 \ds\frac {|u_{\mu_2}|^2}
{M_{\mu_1}}\ud x.
\end{array}\]
We rewrite the last term as follows 
\[
\ds\int_{\mathbb R^N}|\nabla\Phi|^2\ds\frac{M_{\mu_2}}{M_{\mu_1}} \ds\frac {|u_{\mu_2}|^2}
{M_{\mu_2}}\ud x,
\]
where
\[|\nabla\Phi|^2\ds\frac{M_{\mu_2}}{M_{\mu_1}}=|\nabla\Phi|^2e^{(\mu_1-\mu_2)\Phi}
\]
lies in $L^\infty(\mathbb R^N)$ since $\mu_2>\mu_1$. 
 Indeed, we use \eqref{hPhi2}. On the one hand, we write
$\nabla\Phi(x)=\nabla\Phi(0)+\int_0^1 D^2\Phi(\theta x)x\ud \theta$ and, since $D^2\Phi\in L^\infty(\mathbb R^N)$, we get
 $|\nabla\Phi(x)|\leq M|x|$.
On the other hand, $$\Phi(x)=\Phi(0)+\nabla\Phi(0)\cdot x+\int_0^1 (1-\theta)D^2\Phi(\theta x)x\cdot x\ud \theta
\geq \Phi(0)+\ds\frac\Lambda2|x|^2.$$ It follows that 
$|\nabla\Phi(x)|^2e^{(\mu_1-\mu_2)\Phi(x)}\leq M^2|x|^2e^{(\mu_1-\mu_2)(\Phi(0)+\Lambda |x|^2/2)}\leq \frac{C}{\mu_2-\mu_1}$.
We obtain
\[\begin{array}{l}
\ds\frac\gamma2\ds\int_{\mathbb R^N} \ds\frac{\epsilon^2}{M_{\mu_1}}\ud x
+\ds\frac12\ds\int_{\mathbb R^N} M_{\mu_1}\Big|\nabla\Big(\ds\frac{\epsilon^2}{M_{\mu_1}}\Big)\Big|^2\ud x
\\
\leq \ds\frac{(\mu_2-\mu_1)^2}{2\gamma}\ds\int_{\mathbb R^N}\ds\frac{|S|^2}{M_{\mu_1}}\ud x
+\frac{C(\mu_2-\mu_1)}{2\gamma^2}\ds\int_{\mathbb R^N} \ds\frac {|S|^2}
{M_{\mu_2}}\ud x,
\end{array}\]
where $\int_{\mathbb R^N} \frac{|S|^2}{M_{\mu_j}}\ud x\leq \int_{\mathbb R^N} \frac{|S|^2}{M_{\mu_*}}\ud x<\infty$.
We conclude by using $$\|\epsilon\|_{L^1(\mathbb R^N)}\leq \|M_{\mu_1}\|_{L^1(\mathbb R^N)}^{1/2}
\left(\ds\int_{\mathbb R^N} \ds\frac{\epsilon^2}{M_{\mu_1}}\ud x\right)^{1/2}$$
that $\lim_{\mu_2\to \mu_1}\mathscr F(\mu_2)=\mathscr F(\mu_1)$.
\QED

\begin{lemma}\label{L1}
%We suppose that $x\cdot\nabla\Phi(x)\geq 0$ holds for any $x\in \mathbb R^N$. 
Suppose \eqref{hPhi2} holds.
 Let  $(1+|x|^2)S\in L^1(\mathbb R^N)$ with $S\geq 0$ and  $S\in L^2(\mathbb R^N, 
\frac{\ud x}{M_\mu})$ for any $\mu>0$.  
Let $\delta:\mathbb R^N\to [0,\infty)$ be a non-negative bounded and continuous function
(In particular   
$\delta\in L^2(\mathbb R^N, 
M_\mu)$ for any $\mu>0$.) We suppose that there exists $\eta, r>0$ such that
$\delta (x)\geq \eta$ on $B(0,r)$.  
%  Let $u_\mu$ be given by Proposition~\ref{prop:LxM}-ii) for right hand side $\mu S$.
  Then, 
$\lim_{\mu\to \infty}\mathscr F(\mu)=+\infty$.
\end{lemma}

\noindent
{\bf Proof.} We bear in mind \eqref{massmu} for the zeroth moment on $u_\mu$.  
Similarly, considering the second moment and using integration by parts, we are led to 
\[
\begin{array}{lll}
\gamma \ds\int_{\mathbb R^N} |x|^2u_\mu\ud x&=&\mu \ds\int_{\mathbb R^N} |x|^2S\ud x
+\ds\int_{\mathbb R^N} |x|^2\nabla\cdot(\nabla u_\mu+\mu u_\mu\nabla \Phi )\ud x
\\
&=&
\mu \ds\int_{\mathbb R^N} |x|^2S\ud x
+2N \ds\int_{\mathbb R^N}u_\mu\ud x-2\mu \ds\int_{\mathbb R^N} x\cdot \nabla \Phi u_\mu \ud x
\\
&\leq &
\mu \ds\int_{\mathbb R^N} \Big(\ds\frac{2N}{\gamma}+|x|^2\Big) S\ud x-2\Lambda \mu  \ds\int_{\mathbb R^N} |x|^2u_\mu\ud x,
\end{array}\]
since \eqref{hPhi2} implies $x\cdot \nabla\Phi(x)=x\cdot (\nabla\Phi(x)-\nabla\Phi(0))\geq \Lambda |x|^2$.
It follows that the second moment is bounded uniformly wrt $\mu>0$ since
\[
\ds\int_{\mathbb R^N} |x|^2u_\mu\ud x\leq \ds\frac{\mu}{\gamma + 2\Lambda\mu}  \ds\int_{\mathbb R^N} \Big(\ds\frac{2N}{\gamma}+|x|^2\Big) S\ud x
\leq \ds\frac{1}{2\Lambda}  \ds\int_{\mathbb R^N} \Big(\ds\frac{2N}{\gamma}+|x|^2\Big) S\ud x.\]
We now split, for $r>0$, 
\[\begin{array}{lll}
\ds\int_{\mathbb R^N} \delta u_\mu\ud x&=&\ds\int_{|x|\leq r} \delta u_\mu\ud x+\ds\int_{|x|>r} \delta u_\mu\ud x
\\ & \geq & 
\eta \ds\int_{|x|\leq r}  u_\mu\ud x=\eta \left(\ds\int_{\mathbb R^N}  u_\mu\ud x- \ds\int_{|x|> r}  u_\mu\ud x\right)
\\
& \geq & 
\eta \mu \ds\int_{\mathbb R^N}  S\ud x- \ds\frac{\eta}{r^2}\ds\int_{\mathbb R^N}  |x|^2 u_\mu\ud x
\\&\geq& 
\eta \mu \ds\int_{\mathbb R^N}  S\ud x- \ds\frac{\eta}{2\Lambda r^2}\ds\int_{\mathbb R^N}  \Big(\ds\frac{2N}{\gamma}+|x|^2\Big) S\ud x
\end{array}\]
where the RHS tends to $+\infty$ as $\mu\to \infty$.
\QED

In fact, we can make the behavior for large $\mu$'s more precise, by appealing to the 
Laplace method, which can be summarized in the following claim 
\cite[Theorem~15.2.2]{BaSi}.

\begin{lemma}
Let $f:\mathbb R^N\to \mathbb R$ be a continuous function such that $f(0)\neq 0$.
Then, as $\mu$ goes to $+\infty$, 
$\int_{\mathbb R^N} f M_{\mu}\ud x$ is equivalent to 
\[\ds\frac{f(0)}{\mu^{N/2}}\sqrt{
\ds\frac{(2\pi)^N}{\mathrm{det}(\mathrm D^2\Phi(0))}}, 
\]\end{lemma}

\begin{coro}In particular, assuming $\delta\in L^\infty(\mathbb R^N)$, $\delta $ continuous with $\delta (0)\neq 0$, we have $\mathscr F(\mu)\sim \mu \delta (0)\frac{\langle S\rangle}{\gamma}$ as $\mu\to \infty$.
\end{coro}

\noindent
{\bf Proof.}
Let us set 
$$
m(\mu)=\ds\int_{\mathbb R^N} M_\mu\ud x, \qquad
\varsigma(\mu)=\ds\int_{\mathbb R^N} \ds\frac{S^2}{M_\mu}\ud x$$
and introduce the following rescaling
$$\tilde u_\mu(x)=\ds\frac{u_\mu(x)}{\mu\sqrt{\varsigma(\mu)}}.$$
The latter satisfies
\[\gamma \tilde u_\mu-\nabla\cdot\left(M_\mu\nabla \ds\frac{\tilde u_\mu}{M_\mu}\right)=\ds\frac{S}{\sqrt{\varsigma(\mu)}}.\]
It follows that, by using the elementary inequality $2ab\leq a^2+b^2$,  
\[\begin{array}{lll}
\gamma \ds\int_{\mathbb R^N}\ds\frac{|\tilde u_\mu|^2}{M_\mu}\ud x
+\ds\int_{\mathbb R^N} M_\mu 
\left|\nabla \ds\frac{\tilde u_\mu}{M_\mu}\right|^2\ud x&=&
\ds\int_{\mathbb R^N}\ds\frac{S}{\sqrt{\varsigma(\mu)}}\ds\frac{\tilde u_\mu}{M_\mu}\ud x
\\&\leq& \ds\frac{\gamma}{2} 
 \ds\int_{\mathbb R^N}\ds\frac{|\tilde u_\mu|^2}{M_\mu}\ud x
 +\ds\frac{1}{2\gamma} \ds\int_{\mathbb R^N}\ds\frac{S^2}{\varsigma(\mu) M_\mu}\ud x
.\end{array}\]
In turn, we obtain the following estimate
\[\gamma \ds\int_{\mathbb R^N} \ds\frac{|\tilde u_\mu|^2}{M_\mu}\ud x
+2\ds\int_{\mathbb R^N} M_\mu \left|\nabla \ds\frac{\tilde u_\mu}{M_\mu}\right|^2 \ud x
\leq \ds\frac1\gamma,\]
together with 
\[\langle \tilde u_\mu \rangle = \ds\frac{\langle S\rangle}{\gamma\sqrt{\varsigma(\mu)}}.\]
Owing to \eqref{coer_L}, it leads to 
\[
\ds\int_{\mathbb R^N} \left |\tilde u_\mu-\ds\frac{\langle S\rangle}{\gamma\sqrt{\varsigma(\mu)}} \ds\frac{M_\mu}{m(\mu)}\right|\ud x
\leq \ds\frac1{2\sqrt{\gamma\Lambda\mu}}\xrightarrow [\mu\to \infty]{}0.
\]
Therefore, assuming $\delta\in L^\infty(\mathbb R^N)$, $\delta $ continuous with $\delta (0)\neq 0$, we can evaluate 
\[\begin{array}{lll}
\mathscr F(\mu)&=&\mu \sqrt{\varsigma(\mu)}
\ds\int_{\mathbb R^N} \tilde u_\mu \delta \ud x
\underset{\mu \to \infty}{ \sim}
 \mu \sqrt{\varsigma(\mu)}
\ds\int_{\mathbb R^N} \ds\frac{\langle S\rangle}{\gamma\sqrt{\varsigma(\mu)}} \ds\frac{M_\mu}{m(\mu)} \delta \ud x
\\&\underset{\mu \to \infty}{ \sim}&
\ds\frac{\langle S\rangle}{\gamma} \ds\frac{\mu}{m(\mu)} \ds\frac{\delta (0)}{\mu^{N/2}}\sqrt{\ds\frac{(2\pi)^N}{\mathrm{det}(\mathrm D^2\Phi(0))}}
\underset{\mu \to \infty}{ \sim}
\mu \delta (0)\ds\frac{\langle S\rangle}{\gamma}.\end{array} 
\]\QED
Since the function $\mu\mapsto \mathscr F(\mu)$ is continuous, with $\mathscr F(0)=0$, 
we deduce the following existence result.

\begin{theo}\label{C:ex} Let the assumptions of Lemma~\ref{L1} be fulfilled.
%Suppose  $\delta\in L^\infty(\mathbb R^N)$, $\delta $ continuous with $\delta (0)\neq 0$.
For any $\ell\geq 0$, there exists at least a $\mu\in [0,\infty)$ such that \eqref{def_FP}-\eqref{constraint} holds.
\end{theo}

\noindent 
For analysing further the problem we will use the derivability of $\mathscr F$, 
that already appeared for establishing Theorem~\ref{theo}.

\begin{lemma}\label{lem:C1}
The function $\mu\mapsto F(\mu)$ is derivable with 
$\mathscr F'(\mu)=\int_{\mathbb R^N} \delta u'_\mu\ud x$ where $u'_\mu\in L^1(\mathbb R^N)$ satisfies \[
\gamma u_\mu'-\Delta  u_\mu'-\mu\nabla\cdot (u_\mu'\nabla\Phi)=
S
+\nabla\cdot (u_{\mu}\nabla\Phi).
\]
\end{lemma}

\noindent
{\bf Proof.}
Let us consider the equation that defines $u'_\mu$
\[
\gamma v-\Delta  v-\mu\nabla\cdot (v\nabla\Phi)=
S
+\nabla\cdot (u_{\mu}\nabla\Phi)
\]
which is of course recast as
\[(\gamma-L_\mu)u'_\mu=S
+\nabla\cdot (u_{\mu}\nabla\Phi)
.\]
We introduce the weight $$\omega(x)=e^{-\alpha|x|},$$
and the functional space
\[
H_\omega=\Big\{u\in  L^2\Big(\mathbb R^N,\frac{\omega\ud x}{M_\mu}\Big),\ \nabla\Big( \ds\frac u{M_\mu}\Big)\in L^2(\mathbb R^N,\omega M_\mu \ud x)\Big\}.\]
We first show that $\alpha>0$ can be chosen sufficiently small so that 
 $(\gamma-L_\mu)$ is coercive on $H_\omega$.
 Indeed, we have
 \[
 \ds\int_{\mathbb R^N} (\gamma -L_\mu) u \ds\frac{u}{M_\mu}\omega\ud x
 =
 \gamma  \ds\int_{\mathbb R^N}  \ds\frac{|u|^2}{M_\mu}\omega\ud x
 +
 \ds\int_{\mathbb R^N}  \Big|\nabla \Big(\ds\frac{u}{M_\mu}\Big)\Big|^2 M_\mu \omega\ud x
+
 \ds\int_{\mathbb R^N} M_\mu   \nabla \Big(\ds\frac{u}{M_\mu}\Big)\cdot \ds\frac{u}{M_\mu} \nabla \omega\ud x.\]
 The last term recasts as
 \[
  \ds\int_{\mathbb R^N} \sqrt {\omega M_\mu}   \nabla \Big(\ds\frac{u}{M_\mu}\Big)\cdot \ds\frac{\sqrt\omega u}{\sqrt{M_\mu}} \ds\frac{\nabla \omega}{\omega}\ud x.\]
However $\frac{|\nabla\omega|}{\omega}=\alpha$ and by choosing $0<\alpha\ll 1$ sufficiently   
small, we can find $c_*>0$ such that 
\[
\ds\int_{\mathbb R^N} (\gamma -L_\mu) u \ds\frac{u}{M_\mu}\omega\ud x\geq c_\star
 \left(\ds\int_{\mathbb R^N}  \ds\frac{|u|^2}{M_\mu}\omega\ud x
 +
 \ds\int_{\mathbb R^N}  \Big|\nabla \Big(\ds\frac{u}{M_\mu}\Big)\Big|^2 M_\mu \omega\ud x\right).\]
 Note that how much $\alpha$ should be small does not depend on $\mu$.

 The second step is to justify that $S+\nabla\cdot(u_\mu\nabla \Phi ) $ lies in the dual of $H_\omega$.
 For the source term $S$, this simply relies on the fact that 
$\int_{\mathbb R^N} \frac{|S|^2}{M_\mu}\omega\ud x<\infty$.
Next, we use the following estimates
\[\begin{array}{l}
\left|\ds\int_{\mathbb R^N} \nabla\Phi u_\mu\cdot\omega\nabla\Big(\ds\frac{v}{M_\mu} \Big)\ud x\right|
\leq \|\sqrt\omega\nabla\Phi\|_{L^\infty(\mathbb R^N)}
 \ds\int_{\mathbb R^N}  \ds\frac{|u_\mu|}{\sqrt{M_\mu}} \sqrt{\omega M_\mu} \Big|\nabla\Big(\ds\frac{v}{M_\mu}\Big)\Big| \ud x
\\
\hspace*{3cm} \leq\ds\frac12\|\sqrt\omega\nabla\Phi\|_{L^\infty(\mathbb R^N)}^2  \ds\int_{\mathbb R^N}  \ds\frac{|u_\mu|^2}{M_\mu} \ud x
+\ds\frac12
\ds\int_{\mathbb R^N}\omega M_\mu \Big|\nabla\ds\frac{v}{M_\mu}\Big|^2 \ud x
,
\end{array}\]
and 
\[\begin{array}{l}
\left|\ds\int_{\mathbb R^N} \nabla\Phi u_\mu\cdot\ds\frac{v}{M_\mu}\nabla\omega \ud x\right|
\leq
\ds\int_{\mathbb R^N} \sqrt\omega|\nabla\Phi| 
\ds\frac{|u_\mu|}{\sqrt{M_\mu}}\ds\frac{\sqrt\omega |v|}{\sqrt{M_\mu}}\ds\frac{|\nabla\omega|}{\omega} \ud x
\\
\hspace*{3cm}\leq
\ds\frac12  \ \Big\|\ds\frac{\nabla\omega}{\omega}\Big\|_{L^\infty(\mathbb R^N)}
\left(\|\sqrt\omega\nabla\Phi\|_{L^\infty(\mathbb R^N)}^2\ds\int_{\mathbb R^N}\ds\frac{|u_\mu|^2}{M_\mu}\ud x
+
\ds\int_{\mathbb R^N}\omega
\ds\frac{ |v|^2}{M_\mu} \ud x\right).
\end{array}\]
As already observed $|\nabla\Phi(x)|\leq M|x|$ and thus $\sqrt\omega\nabla\Phi$ lies in $L^\infty(\mathbb R^N)$.
We conclude that $u'_\mu$ is well-defined in $H_\omega$.

Finally, let $h>0$ and set 
$$U_h=\ds\frac{u_{\mu+h}-u_\mu}{h}-u'_\mu.$$
It satisfies 
\[
(\gamma-L_\mu)U_h = \nabla\cdot((u_{\mu+h}-u_\mu)\nabla \Phi).\]
With the same manipulations, we control
\[\ds\int_{\mathbb R^N} 
\omega \ds\frac{|U_h|^2}{M_\mu}\ud x \text{ and }
\ds\int_{\mathbb R^N} 
\omega M_\mu \Big|\nabla\Big( \ds\frac{U_h}{M_\mu}\Big)\Big|^2\ud x
\]
by 
\[
\ds\int_{\mathbb R^N}\ds\frac{|u_{\mu+h}-u_\mu|^2}{M_\mu}\ud x \xrightarrow [h\to 0]{}0.
\]
It implies that $$\lim_{h\to 0}\left(
\ds\frac{\mathscr F(\mu+h)-\mathscr F(\mu)}{h}-\int_{\mathbb R^N} u'_\mu\delta \ud x\right)=
\lim_{h\to 0}\ds\int_{\mathbb R^N} U_h\delta \ud x=0.$$
\QED

\section{Monotonicity}
\label{Sec:Monot}

We are going to show that, in the radial symmetry case and under the compatibility conditions stated in \eqref{cx},
the function $\mathscr F$ is increasing; to this end we use a duality argument. We introduce the solution $\psi_\mu$ of 
\[
(\gamma-L^\star_\mu) \psi_\mu= \delta\]
 so that  $\mathscr F(\mu) 
$ recasts as \[
\ds\int_{\mathbb R^N} u_\mu\delta\ud x 
=
\ds\int_{\mathbb R^N} u_\mu(\gamma-L^\star_\mu) \psi_\mu\ud x 
=\ds\int_{\mathbb R^N} (\gamma-L_\mu) u_\mu\psi_\mu\ud x 
=\mu\ds\int_{\mathbb R^N} S\psi_\mu\ud x.
\]
Accordingly, showing the monotonicity of  $\mathscr F$
reduces to investigating the sign of 
\begin{equation}\label{deriv}
\ds\frac{\ud}{\ud \mu} \mathscr F(\mu) =\ds\int_{\mathbb R^N} S\psi_\mu\ud x+\mu\ds\int_{\mathbb R^N} S\psi'_\mu\ud x\end{equation}
where $\psi'_\mu$ satisfies
\[(\gamma-L^\star_\mu) \psi'_\mu=-\nabla\Phi\cdot \nabla\psi_\mu.\]
The data $S$ and $\delta$ being non-negative, the first integral in the right hand side of \eqref{deriv} is non-negative. We are going 
to show that the second term is equally non-negative, under appropriate assumptions on the data. 

\subsection{Problem in radial symmetry}

From now on, we assume  that all data $S,\delta, \Phi$ are  \emph{radially symmetric}.
In turn, $M_\mu$ and the  solutions of the associated PDE are also radially symmetric.
We write the equation in radial coordinates: we get 
\[
\gamma u_\mu -\mu \partial_r (u_\mu\partial_ r \Phi)-\mu\ds\frac{N-1}{r}\partial_r \Phi u_\mu-\ds\frac{1}{r^{N-1}}\partial_r(r^{N-1}\partial_r u_\mu)=\mu S.\]
It casts as 
\[
\gamma u_\mu -\ds\frac{1}{r^{N-1}}\partial_r \left(M_\mu r^{N-1} \partial_r\Big(\ds\frac{u_\mu}{M_\mu}\Big)\right)
=\mu S.\]
For the adjoint equation, we obtain
\begin{equation}\label{adjrad}\begin{array}{l}
\gamma \psi_\mu +\mu \partial_r \Phi\partial_r \psi_\mu -\ds\frac{1}{r^{N-1}}\partial_r(r^{N-1}\partial_ r \psi_\mu)=\delta
\\
\qquad\qquad
= \gamma \psi_\mu - \ds\frac{1}{r^{N-1} M_\mu}\partial_r(r^{N-1}M_\mu\partial_r \psi_\mu).
\end{array}\end{equation}
Let us set $\chi_\mu=\partial_r\psi_\mu$.
It satisfies
\[\begin{array}{l}
\left(\gamma +
\ds\frac{N-1}{r^2}+\mu \partial^2_r \Phi\right)\chi_\mu +\mu \partial_r \Phi\partial_r \chi_\mu 
-\ds\frac{1}{r^{N-1}}\partial_r(r^{N-1}\partial_ r \chi_\mu)=\partial_r\delta
\\\qquad\qquad
=\left(\gamma +
\ds\frac{N-1}{r^2}+\mu \partial^2_r \Phi\right)\chi_\mu  
-\ds\frac{1}{r^{N-1}M_\mu }\partial_r(r^{N-1}M_\mu \partial_ r \chi_\mu).\end{array}\]
From now on, we  assume the convexity/monotonicity properties \eqref{cx}.
%\begin{equation}\label{cx}\partial_r \Phi\geq 0,\qquad \partial^2_r \Phi\geq 0,\qquad \partial_r\delta\leq 0.\end{equation}
Under these assumptions, by the maximum principle, we obtain $$\chi_\mu=\partial_r\psi_\mu \leq 0.$$

Now, we go back to \eqref{adjrad}, which yields
\[
\begin{array}{l} \gamma\psi'_\mu +\mu \partial_r \Phi\partial_r \psi'_\mu -\ds\frac{1}{r^{N-1}}\partial_r(r^{N-1}\partial_ r \psi'_\mu)=- 
\underbrace{\partial_r \Phi}_{\geq 0}
\underbrace{\partial_r \psi_\mu }_{\leq 0}
\\
\qquad\qquad
= \gamma \psi'_\mu - \ds\frac{1}{r^{N-1} M_\mu}\partial_r(r^{N-1}M_\mu\partial_r \psi'_\mu).\end{array}\]
The right hand side thus satisfies 
$-\partial_r \Phi\partial_r \psi_\mu\geq 0$ so that $ \psi'_\mu\geq 0$.
Coming back to \eqref{deriv}, we conclude that $\mathscr F$ is non decreasing.

We need to slightly improve the result, requiring further regularity on $\delta, \Phi$, say 
$\delta\in C^{1}$, $\Phi\in C^2$, with $\delta, \Phi$ not identically 0. We can thus apply the strong maximum principle \cite[Section~3.2]{GT}
which tells us that $\psi'_\mu>0$ on $(0,\infty)$.
Our findings recap as follows.

\begin{theo}\label{mainth}
In the radially symmetric framework, we suppose that $S, \delta\in C^1$ and $\Phi\in C^2$  take non-negative values,
but are not identically 0,
and  that \eqref{cx} is fulfilled.
Then, $\mathscr F$ is increasing and, for any $\ell\geq 0$, the problem  \eqref{def_FP}-\eqref{constraint}  admits a unique solution $0<\mu <\infty$.
\end{theo}

The assumptions of radial symmetry together with \eqref{cx} are quite natural and relevant  for the presented modeling: the tumor being located at $x=0$, the action of the immune cells, embodied into $\delta$, is centred  on this position and the chemotactic potential $\Phi$ 
 drives the immune cells towards the tumoral centre. Let us detail a simple example showing that these assumptions are also technically important.
 We consider the simplest confining potential $\Phi(x)=|x|^2$ and we compute the first even moments 
of the solutions of \eqref{def_FP}:
we have already seen that $\gamma\int_{\mathbb R^N} u_\mu\ud x=\mu \int_{\mathbb R^N} S\ud x$; next we have
\[\begin{array}{lll}
\gamma\ds\int_{\mathbb R^N} |x|^2 u_\mu\ud x&=& \mu \ds\int_{\mathbb R^N} |x|^2 S\ud x
+2N \ds\int_{\mathbb R^N} u_\mu \ud x- 2\mu \ds\int_{\mathbb R^N} u_\mu x\cdot\nabla\Phi\ud x\\
&=& \mu \ds\int_{\mathbb R^N} |x|^2S\ud x
+2N\ds\frac\mu\gamma  \ds\int_{\mathbb R^N} S \ud x-4 \mu \ds\int_{\mathbb R^N}  |x|^2 u_\mu\ud x
\end{array}\]
so that 
\[
\ds\int_{\mathbb R^N} |x|^2 u_\mu\ud x=\ds\frac{\mu}{\gamma + 4\mu}
 \ds\int_{\mathbb R^N} (|x|^2+2N/\gamma) S\ud x.
\]
It follows that $\mu\mapsto \int_{\mathbb R^N} |x|^2 u_\mu\ud x$ is an increasing function.
We turn to 
\[
\gamma\ds\int_{\mathbb R^N} |x|^4 u_\mu\ud x= \mu \ds\int_{\mathbb R^N} |x|^4 S\ud x
+(4N+8) \ds\int_{\mathbb R^N} |x|^2 u_\mu \ud x- 8\mu \ds\int_{\mathbb R^N} |x|^4 u_\mu \ud x.
\]
It leads to the expression 
\[
\ds\int_{\mathbb R^N} |x|^4 u_\mu\ud x=\underbrace{A\frac{\mu}{\gamma + 8\mu}+ B\frac{1}{\gamma + 8\mu}\frac{\mu}{\gamma + 4\mu}}_{:=f(\mu)}\]
with 
\[A=\ds\int_{\mathbb R^N} |x|^4 S\ud x,\qquad
B=(4N+8)\ds\int_{\mathbb R^N} (|x|^2+2N) S\ud x
.\]
The forth momentum is not necessarily a monotone function of $\mu$ since
\[
f'(\mu)=\ds\frac{1}{(\gamma + 8\mu)^2}\left(\gamma A+B\ds\frac{\gamma}{\gamma+4\mu}
-4B\ds\frac{\mu(\gamma+8\mu)}{(\gamma+4\mu)^2}\right)\]
might change sign.
Therefore, if we set 
 $\delta(x)=|x|^4$, 
 the monotonicity of $\mu\mapsto \int_{\mathbb R^N} \delta u_\mu\ud x$ does not hold in general.
 Note however that this example $\delta(x)=|x|^4$
  vanishes at $x=0$, contradicting the modelling assumptions.

\subsection{Numerical illustrations}

Dealing with the radially symmetric problem, the finite elements framework  is a reliable way to get rid of the singularity at $r=0$.  
For realizing simulations, we consider the problem set on the slab $[0,1]$: since the phenomena are naturally quite concentrated next to the origin, we expect that this does not influence too much the final results (by the way, we indeed do not observe
 significant differences 
 when imposing Dirichlet or Neumann conditions at $r=1$ or extending the domain for larger $r$'s).
 We introduce a discretization of $ [0,1]$ with $N+1$ points
\[r_0=0<r_1=h<...<r_N=Nh=1,\qquad h=1/N.\]
We introduce the associated $\mathbb P_1$ basis functions, $\chi_1,...,\chi_{N-1}$:
\[\chi_j(r)=\ds\frac{r-(j-1)h}h\mathbf 1_{(j-1)h<r< jh}+ \ds\frac{(j+1)h-r}h\mathbf 1_{jh\leq r<(j+1)h},\quad
\chi_0(r)=- \ds\frac{r-h}h\mathbf 1_{0\leq r<h}.\]
Then, we define the matrices
with coefficients \[
M_{ij}=\ds\int_0^1\chi'_i(r)\chi'_j(r) r^{n-1}\ud r,\quad 
A_{ij}=\ds\int_0^1\chi_i(r)\chi_j(r) r^{n-1}\ud r.\]
Given the potential $\Phi$, we also define the  centered difference matrix $C$, which is skew-symmetric with
\[
C_{j,j+1}=\ds\frac12((j+1)h)^{n-1}\Phi'((j+1)h).
\]
Then, for a given source term $S$, we define the vector with components
\[S(j)=\ds\int_{ih}^{(i+1)h} S(r)r^{n-1}\ud r.\]
Eventually, we solve the linear system
\[(\gamma A -\mu C+M)U=\mu S,\]
and we compute the associated discrete version of $\mathscr F(\mu)=\int_0^1 \delta u_\mu(r)r^{n-1}\ud r$.

We perform the simulation with a source term given by 
\[S(r)=\mathbf 1_{.3\leq r\leq .5}\]
which, for the application to tumor-immune system interactions,  corresponds to a located 
reservoir of resting immune cells, for instance a blood vessel or a lymph node.
We set $n=3$ and $\gamma=1$.

\begin{figure}[!hbtp]
\begin{center}
\includegraphics[height=6cm]{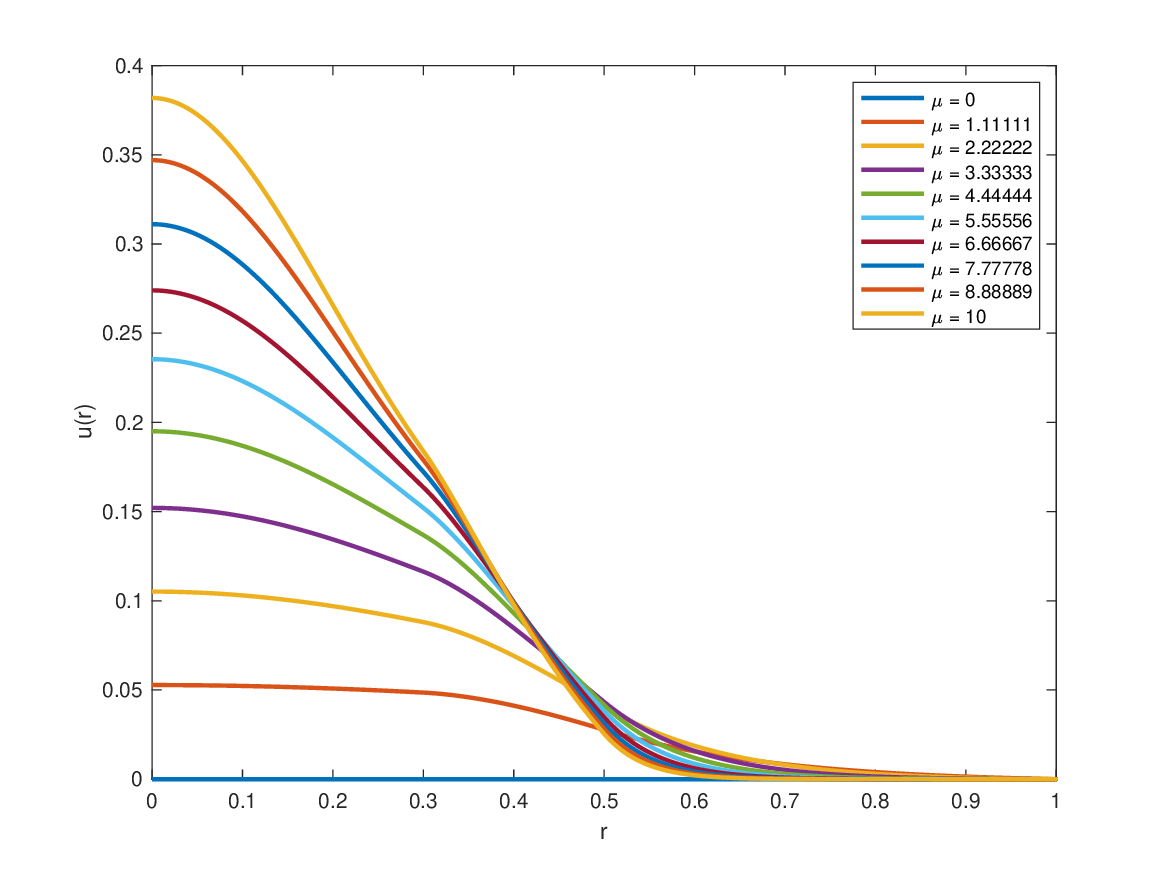}
\caption{Solutions $r\mapsto u_\mu(r)$ for 10 equidistant values of $\mu$ in $[0,10]$. As $\mu$ increases
the solution takes larger value at the origin and presents a stiffer profile for transient radius \label{F_OK1}}
\end{center}
\end{figure}

\begin{figure}[!hbtp]
\begin{center}
\includegraphics[height=6cm]{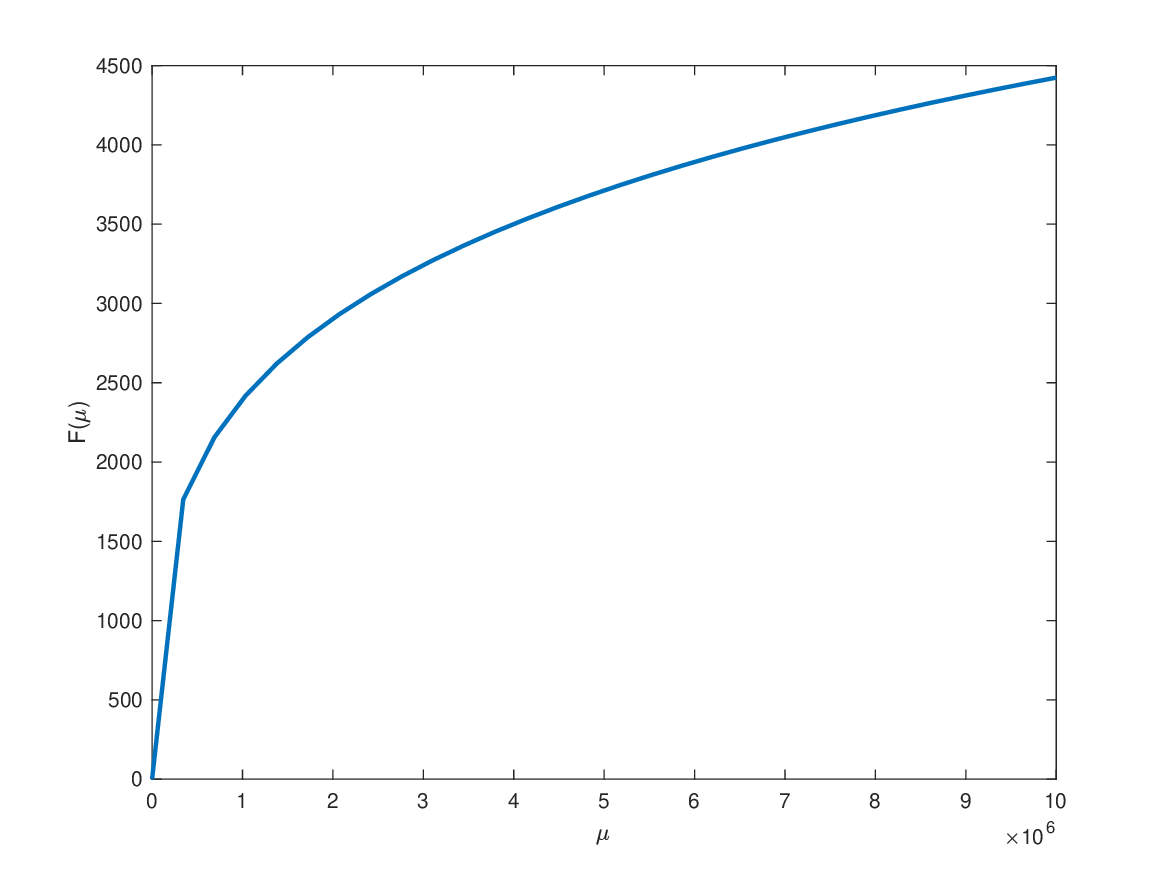}
%~\includegraphics[height=6cm]{./Fig/UMu_OK.eps}
\caption{
%Left: 
$\mu\mapsto \mathscr F(\mu)$ for $\mu$ up to $10^7$.
%Right: Snapshot of the corresponding solution profiles $r\mapsto u_\mu(r)$
\label{F_OK2}}
\end{center}
\end{figure}

We start with simulations of the expected situation, with a confining potential 
pointing towards the origin $\Phi(r)=2r^2$, and 
a constraint kernel 
 peaked at the origin $$\delta(r)=\frac{e^{-r^2/\epsilon}}{(4\pi\epsilon)^{n/2}},\qquad \epsilon=10^{-3}.$$
Fig.~\ref{F_OK1} represents the profile of the solutions for  
relatively small values of $\mu$: as $\mu$ increases, the value at $r=0$ increases
and the solution concentrates near $r=0$.
We numerically check that $\mu\mapsto \mathscr F(\mu)$ is non decreasing, 
see Fig.~\ref{F_OK2}; the solution $u_\mu$ becomes highly concentrated to $r=0$, which thus requires a very fine mesh to resolve the solution when $\mu$ becomes large.

\begin{figure}[!hbtp]
\begin{center}
\includegraphics[height=4cm]{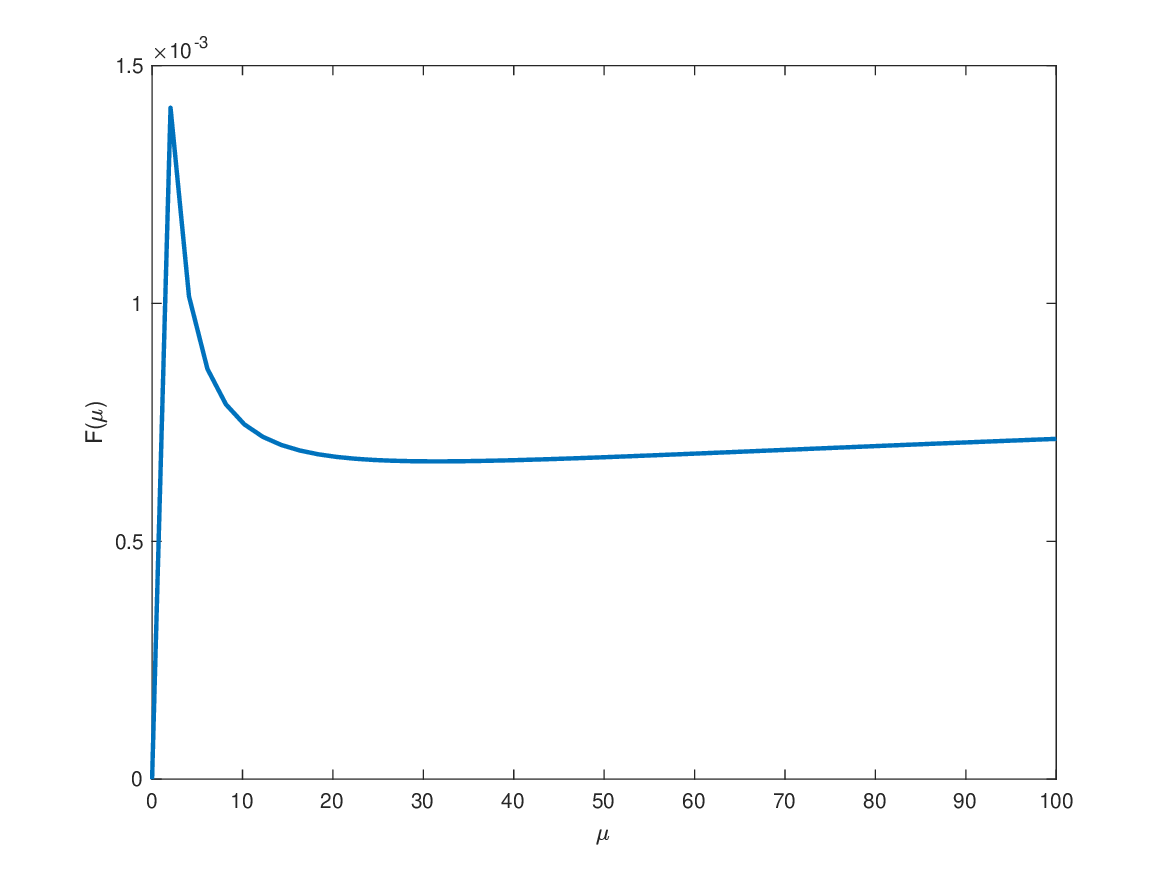}~\includegraphics[height=4cm]{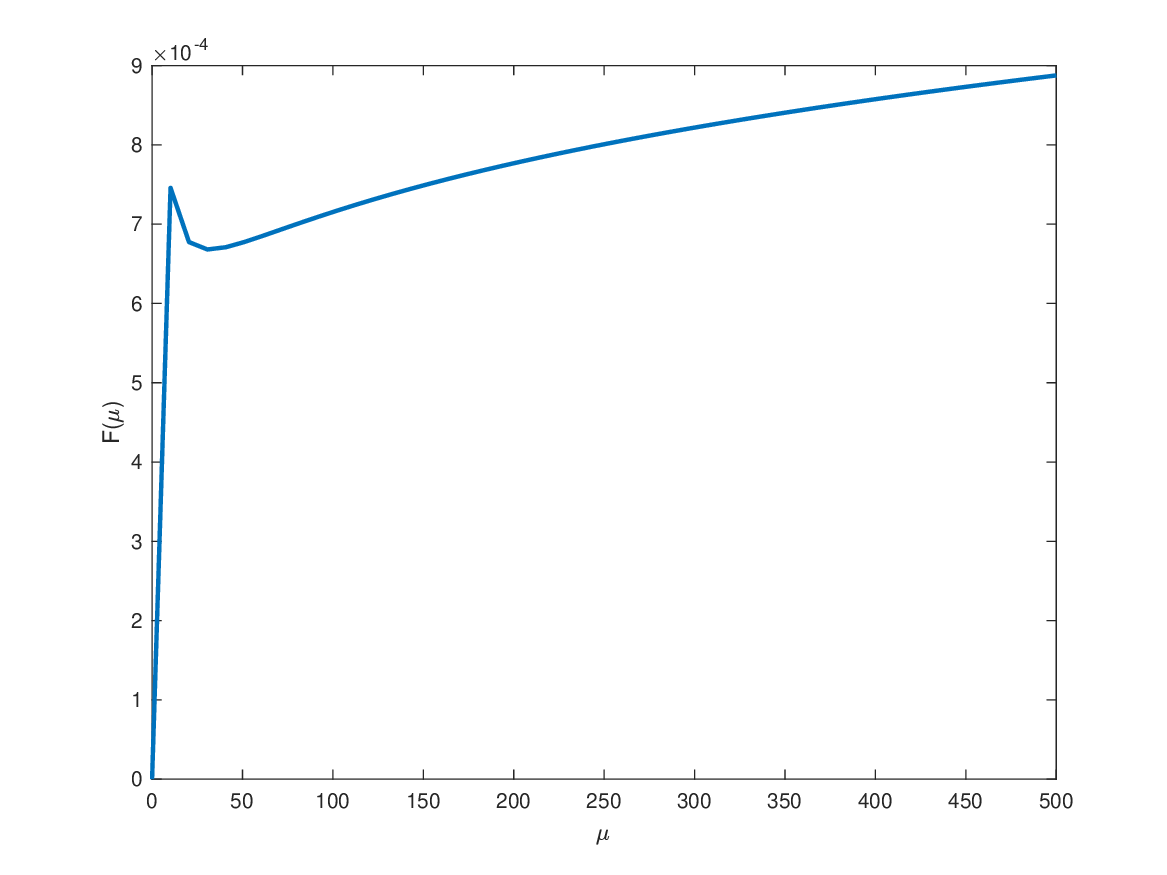}~\includegraphics[height=4cm]{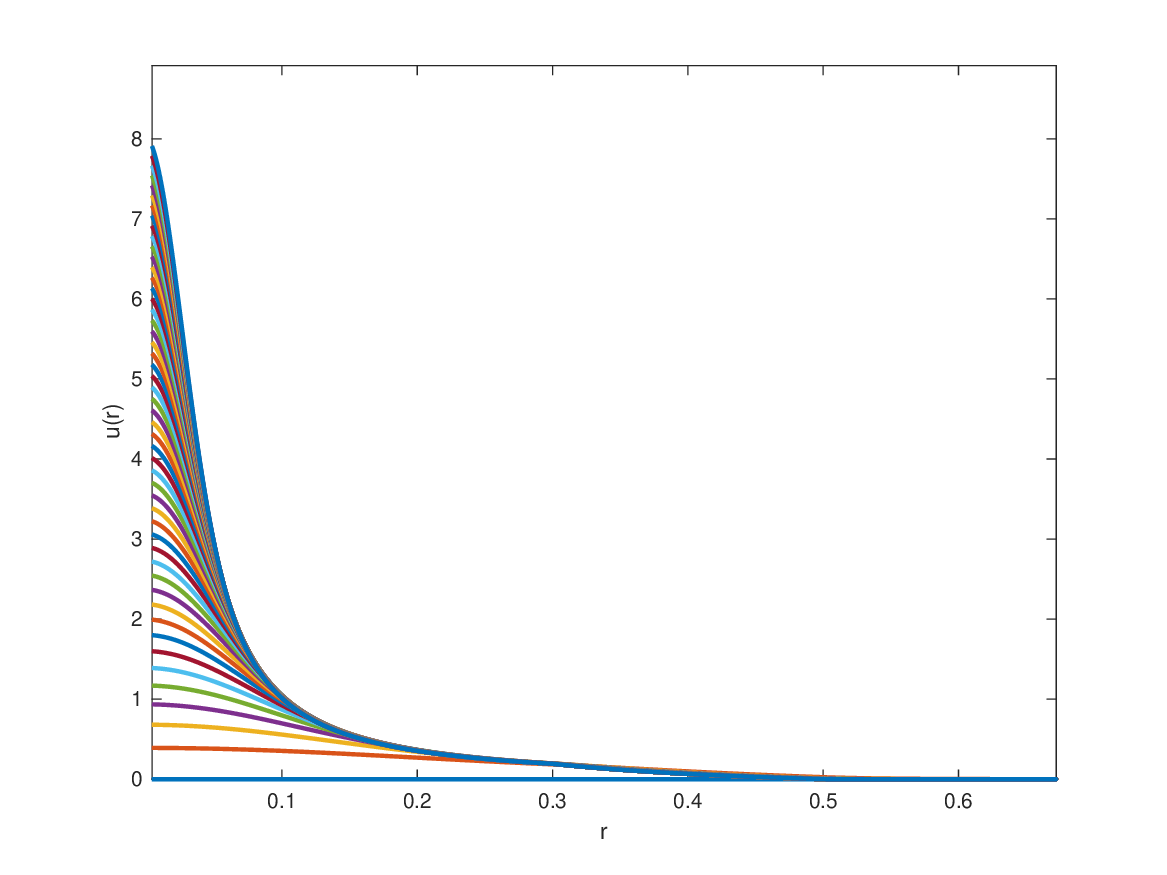}
\caption{Profile of
$\mu\mapsto \mathscr F(\mu)$ for $\mu$ up to $100$ (left),   up to $500$ (middle), 
and snapshots of the corresponding solution profiles $r\mapsto u_\mu(r)$ (right)
\label{F_delta4}}
\end{center}
\end{figure}

Simulations of the
 counter-example detailed in the previous section are displayed in Fig.~\ref{F_delta4}:
 we just modify $\delta$ into
 \[  \delta(r)=\epsilon r^4,\qquad \epsilon=10^{-3}.\]
 Then the function $\mu\mapsto \mathscr F(\mu)$ looses its monotony, but it is still increasing at $\mu=0$ and for large values of $\mu$.
 
\begin{figure}[!hbtp]
\begin{center}
\includegraphics[height=4cm]{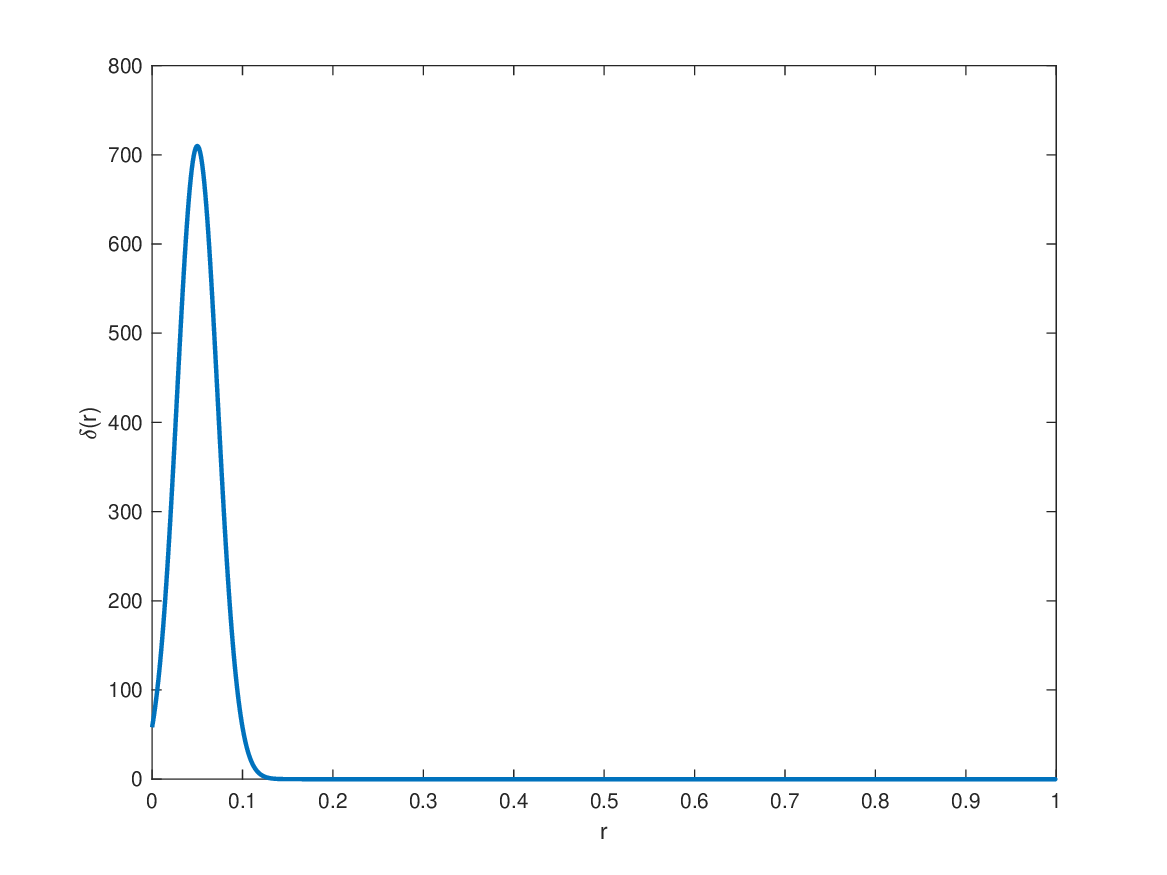}~\includegraphics[height=4cm]{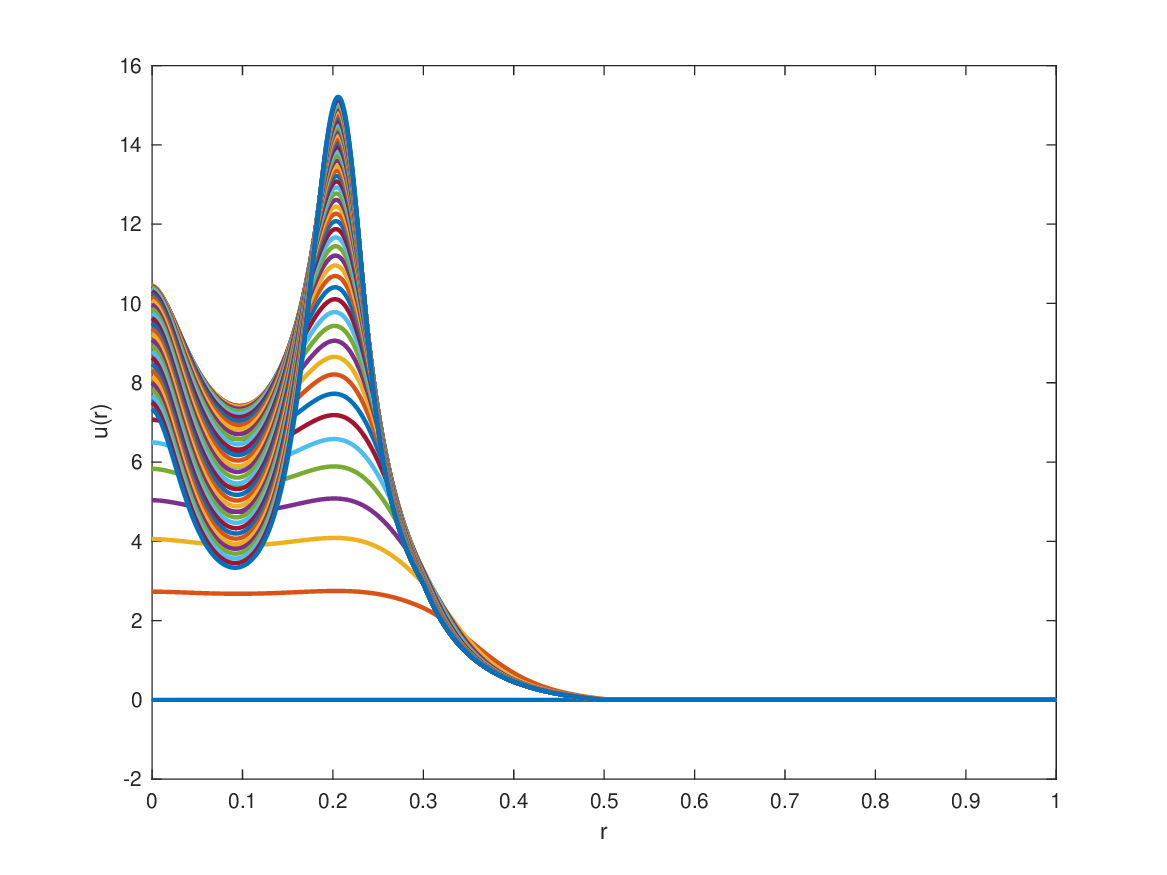}
\\
\includegraphics[height=4cm]{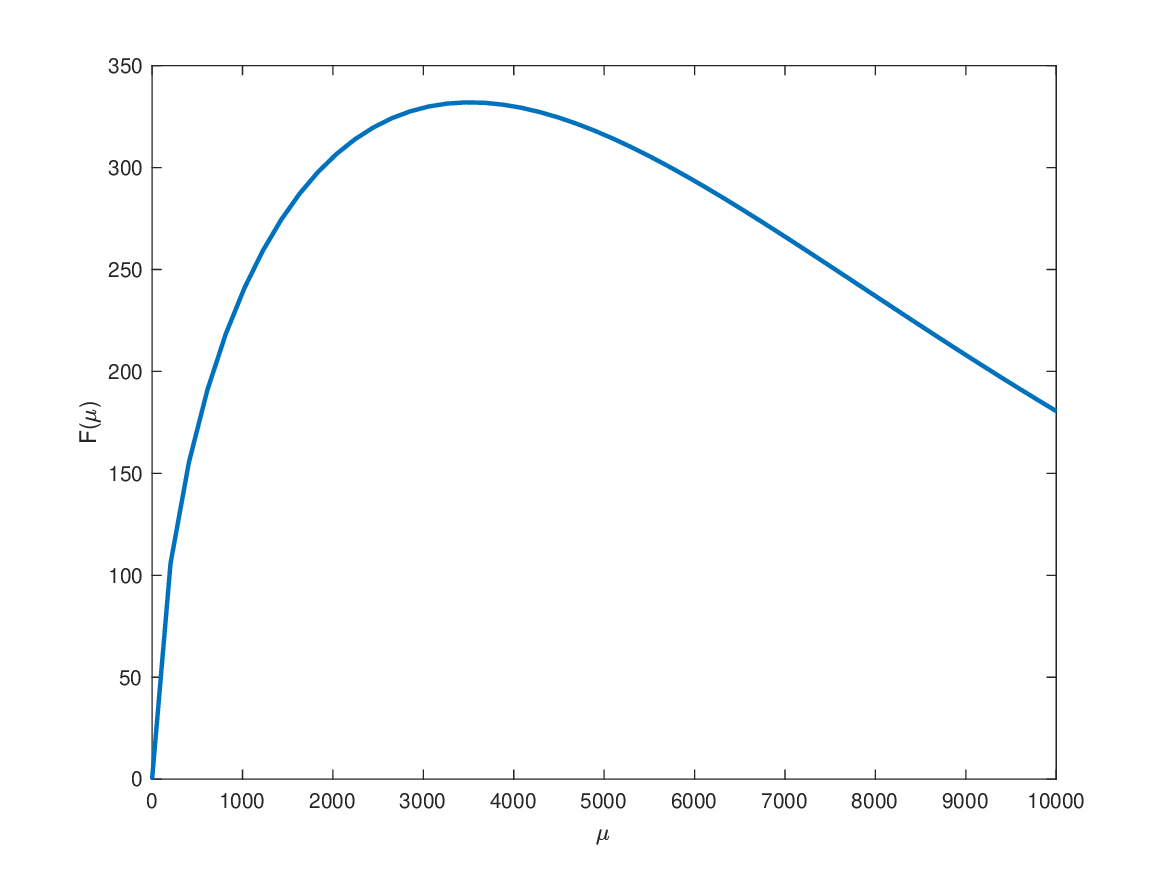}~\includegraphics[height=4cm]{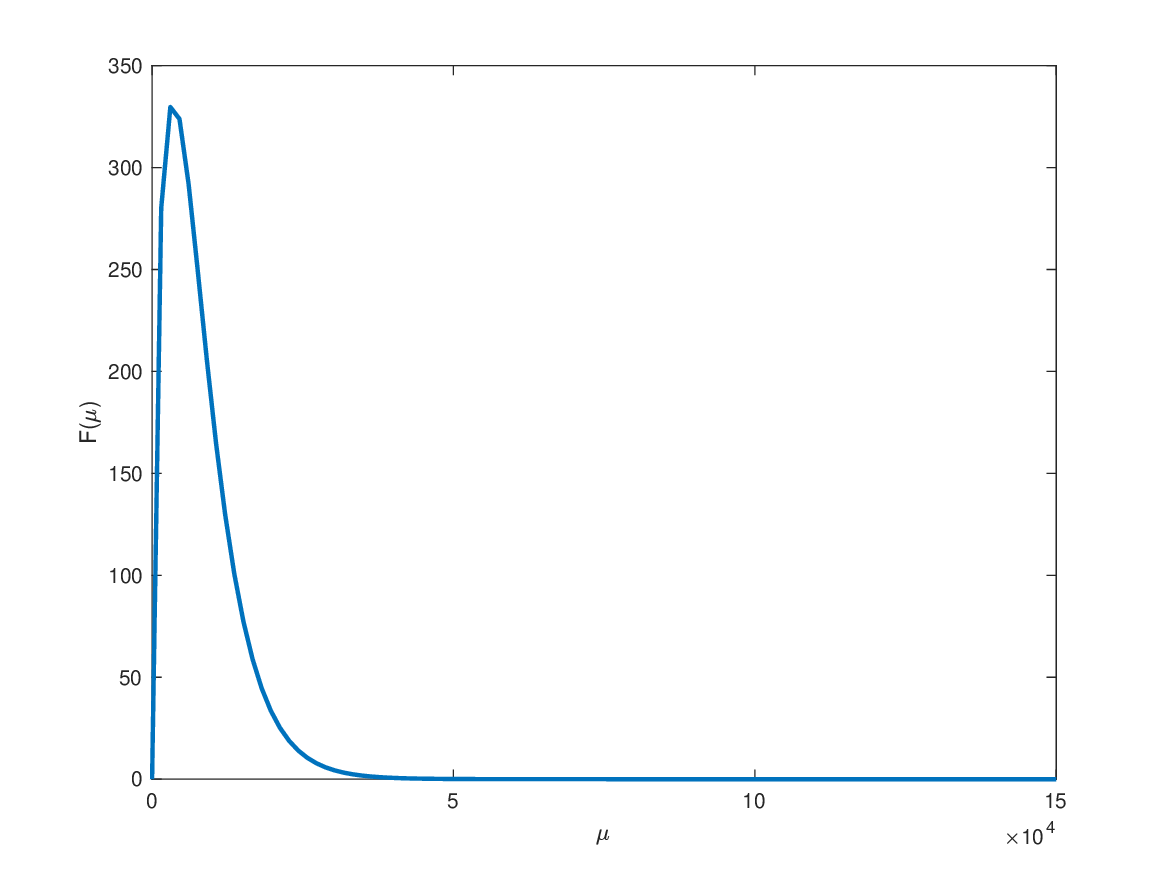}
\caption{Profile of
$r\mapsto \delta(r)$ (top-left), 
profile of $\mu\mapsto \mathscr F(\mu)$ for $\mu$ up  to $10^4$ (bottom-left) and $15\cdot 10^4$
(bottom-right), 
the  solution profiles $r\mapsto u_\mu(r)$ for several $\mu$ up  to $10^4$ (top-right) \label{F_Bad}}
\end{center}
\end{figure}

Next, we consider a situation which can find some physical motivation: 
the potential is given by 
\[
\Phi(r)=2r^2\times (r-r_c)^2,\qquad r_c=.2.\]
It has a quite flat profile between the two minima $r=0$ and $r=.2$.
Note that 
\eqref{cx} is not satisfied: $\Phi$ is not convex, and not monotone.
It describes a defect of the  attractivity of the immune cells towards the tumor, due either 
to the geometry of the tissues around the tumor, or to pro-tumoral effects that reduce the efficacy of the immune response.
Another pro-tumoral effect can result in a reduction of the capacity of the immune cells to eliminate tumor cells, 
that we traduce by shifting the kernel $\delta$
\[\delta(r) =
\frac{e^{-(r-r_1)^2/\epsilon}}{(4\pi\epsilon)^{n/2}},\qquad \epsilon=10^{-3},\qquad r_1=.05.\]
Note that it still keeps a significantly  positive 
value at $r=0$, see Fig.~\ref{F_Bad}-Top Left.
We indeed observe 
that $\mu\mapsto \mathscr F(\mu)$ does not tend to infinity as $\mu \to \infty$, 
and the monotonicity is compromised, see Fig.~\ref{F_Bad}-Bottom.
The solution $u_\mu$ tends to form a high peak in the interior of the domain, thus far from the tumor, Fig.~\ref{F_Bad}-Top Right.

\begin{figure}[!hbtp]
\begin{center}
\includegraphics[height=4cm]{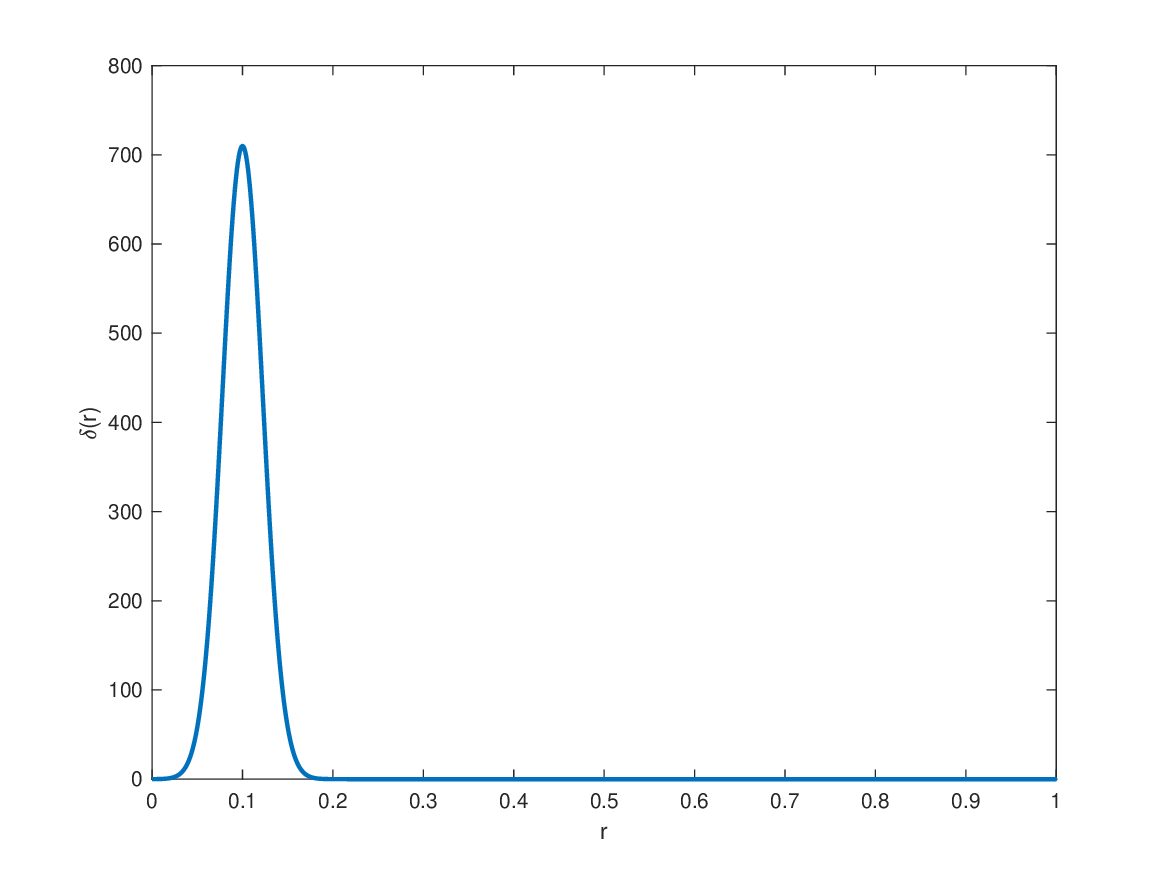}~\includegraphics[height=4cm]{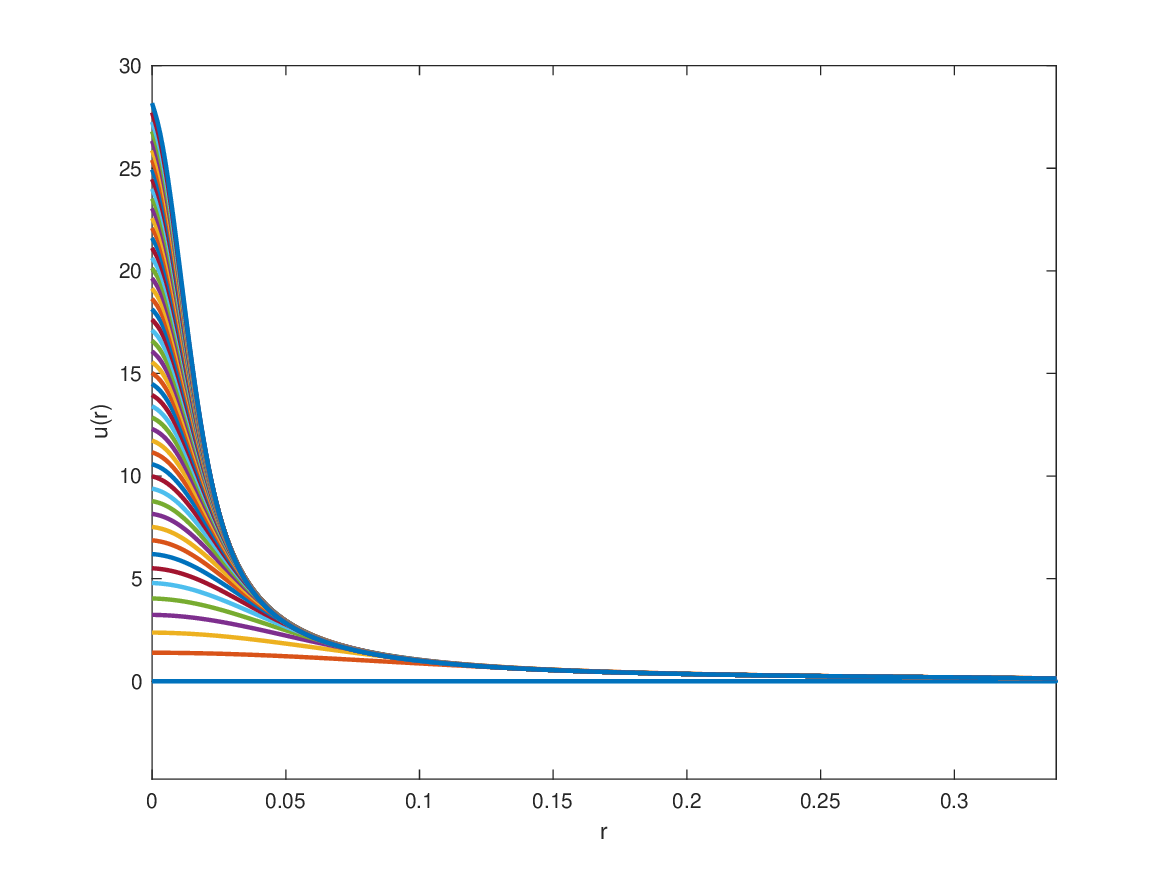}~
\includegraphics[height=4cm]{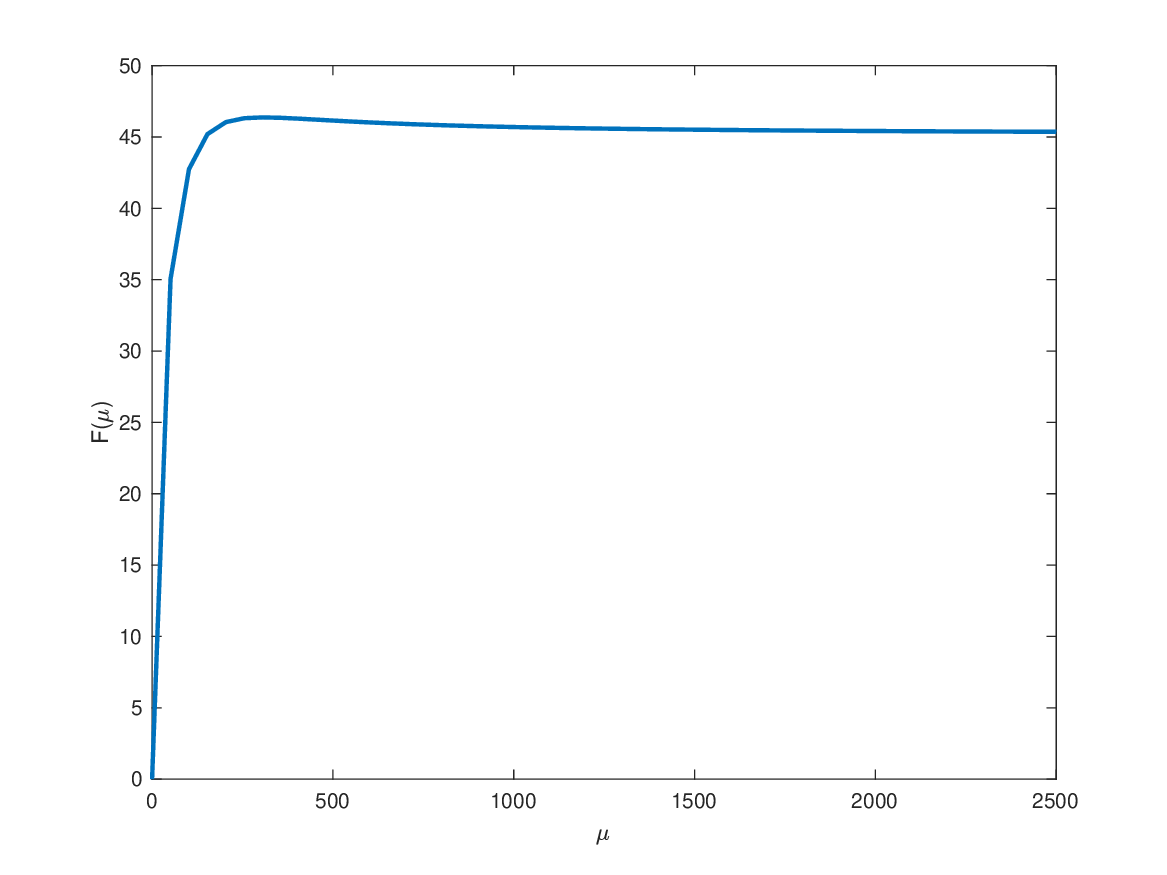}
\caption{Profile of
$r\mapsto \delta(r)$ (left), 
profile of $\mu\mapsto \mathscr F(\mu)$ for $\mu$ up  to $2500$ (bottom-left) 
(bottom-right), snapshot on 
the  corresponding solution profiles $r\mapsto u_\mu(r)$ (top-right)
\label{F_Bad2}}
\end{center}
\end{figure}

Eventually, we challenge the condition that $\delta$ takes positive value near the origine: we come back to the quadratic potential, but now we work with 
\[\delta(r) =
\frac{e^{-(r-r_1)^2/\epsilon}}{(4\pi\epsilon)^{n/2}},\qquad \epsilon=10^{-3},\qquad r_1=.1\]
which (almost) vanishes at $r=0$. Results are 
reported in Fig.~\ref{F_Bad2}. Again, we observe that $\mu\mapsto \mathscr F(\mu)$ is not monotone and does not tend to $\infty$ (it seems to be decaying for large $\mu$'s), at least as far as it can be numerically checked.

These numerical experiments highlight the role of the assumptions on both the potential, and the 
constraint kernel.
Coming back to the motivation from the modeling of tumor-immune system interactions,
these findings shed light on the role of the pro-tumor mechanisms
which not only may promote tumor proliferation, but can also 
reduce the efficacy of the immune response, and eventually allow the tumor to escape to the control of the immune system 
\cite{aabg_pone}.

\section*{Ackowledgements}
This research is supported by the  CNRS program ``International Emerging Actions''.
We acknowledge the support of the   Math.~Dept. at Penn State University.

\bibliography{AGJ}

\begin{thebibliography}{10}

\bibitem{AMTU}
A.~Arnold, P.~Markowich, G.~Toscani, and A.~Unterreiter.
\newblock On convex {S}obolev inequalities and the rate of convergence to
  equilibrium for {F}okker-{P}lanck type equations.
\newblock {\em Comm. Partial Differential Equations}, 26(1-2):43--100, 2001.

\bibitem{aabg_pub}
K.~Atsou, F.~Anju\`ere, V.~M. Braud, and T.~Goudon.
\newblock A size and space structured model describing interactions of tumor
  cells with immune cells reveals cancer persistent equilibrium states in
  tumorigenesis.
\newblock {\em J. Theor. Biol.}, 490:110163, 2020.

\bibitem{aabg_pone}
K.~Atsou, F.~Anju\`ere, V.~M. Braud, and T.~Goudon.
\newblock A size and space structured model of tumor growth describes a key
  role for protumor immune cells in breaking equilibrium states in
  tumorigenesis.
\newblock {\em PlosOne}, page 0259291, 2021.

\bibitem{aabg_fr}
K.~Atsou, F.~Anju\`ere, V.~M. Braud, and T.~Goudon.
\newblock Analysis of the equilibrium phase in immune-controlled tumors
  predicts best strategies for cancer treatment.
\newblock {\em Frontiers in Oncology, Advances in Math. and Comput. Oncology},
  12:878827, 2022.

\bibitem{Bac}
F.~Baccelli, D.R. McDonald, and J.~Reynier.
\newblock A mean field model for multiple {TCP} connections through a buffer
  implementing {RED}.
\newblock {\em Performance Evaluation}, 49(1-4):77--97, 2002.

\bibitem{DGL}
A.~Devys, T.~Goudon, and P.~Lafitte.
\newblock A model describing the growth and size distribution of multiple
  metastatic tumors.
\newblock {\em Disc. Cont. Dyn. Syst.-B}, 12:731--767, 2009.

\bibitem{DoGa}
M.~Doumic-Jauffret and P.~Gabriel.
\newblock Eigenelements of a general aggregation-fragmentation model.
\newblock {\em Math. Models Methods Appl. Sci.}, 20(5):757--783, 2010.

\bibitem{GT}
D.~Gilbarg and N.~Trudinger.
\newblock {\em Elliptic Partial Differential Equations of Second Order}.
\newblock Springer, 1998.

\bibitem{Mich1}
P.~Michel.
\newblock Existence of a solution to the cell division eigenproblem.
\newblock {\em Models Math. Meth. App. Sci.}, 16(Suppl. issue 1):1125--1153,
  2006.

\bibitem{MMP}
P.~Michel, S.~Mischler, and B.~Perthame.
\newblock General relative entropy inequality: an illustration on growth
  models.
\newblock {\em J. Math. Pures et Appl.}, 84(9):1235--1260, 2005.

\bibitem{PerBk}
B.~Perthame.
\newblock {\em Transport equations in biology}.
\newblock Frontiers in Math. Birkhauser, 2007.

\bibitem{PeRy}
B.~Perthame and L.~Ryzhik.
\newblock Exponential decay for the fragmentation or cell-division equation.
\newblock {\em J. Differential Equations}, 210:155--177, 2005.

\bibitem{BaSi}
B.~Simon.
\newblock {\em Advanced Complex Analysis. A Comprehensive Course in Analysis,
  Part 2{B}}.
\newblock AMS, 2015.

\end{thebibliography}
\bibliographystyle{plain}

\end{document}